\documentclass[11pt,reqno,tbtags]{amsart}
\usepackage{amssymb}
\pdfoutput=1
\usepackage[colorlinks=true, pdfstartview=FitV, linkcolor=blue,
  citecolor=blue, urlcolor=blue,pagebackref=false]{hyperref}
\usepackage[numeric,initials,nobysame]{amsrefs}
\usepackage{amsmath,amssymb,amsfonts,dsfont}
\usepackage{enumerate,graphicx}

\numberwithin{equation}{section}

\newtheorem{maintheorem}{Theorem}

\newtheorem{theorem}{Theorem}[section]
\newtheorem*{theorem*}{Theorem}
\newtheorem{lemma}[theorem]{Lemma}
\newtheorem{claim}[theorem]{Claim}
\newtheorem{proposition}[theorem]{Proposition}
\newtheorem{corollary}[theorem]{Corollary}

\newtheorem{conjecture}[theorem]{Conjecture}
\newtheorem{question}[theorem]{Question}

\theoremstyle{definition}{

}

\theoremstyle{remark}{
\newtheorem{remark}[theorem]{Remark}
\newtheorem*{remark*}{Remark}

}

\newcommand{\R}{\mathbb R}

\newcommand{\Z}{\mathbb Z}

\newcommand{\E}{\mathbb{E}}
\renewcommand{\P}{\mathbb{P}}
\DeclareMathOperator{\var}{Var} \DeclareMathOperator{\Cov}{Cov} 

\renewcommand{\epsilon}{\varepsilon}

\newcommand{\given}{\, \big| \,}
\newcommand{\one}{\boldsymbol{1}}

\newcommand{\deq}{\stackrel{\scriptscriptstyle\triangle}{=}}
\newcommand{\tmix}{t_\textsc{mix}}

\newcommand{\gap}{\text{\tt{gap}}}
\newcommand{\sob}{\alpha_{\text{\tt{s}}}}
\newcommand{\lup}{{\ell}}
\newcommand{\ldown}{{r}}
\newcommand{\allplus}{{\underline{1}}}
\DeclareMathOperator{\arctanh}{arctanh}
\DeclareMathOperator{\dist}{dist}

\begin{document}

\title[Mixing of critical Ising model on trees]{Mixing time of critical Ising model on trees \\ is polynomial in the height}
\date{}
\author{Jian Ding, \thinspace Eyal Lubetzky and Yuval Peres}

\address{Jian Ding\hfill\break
Department of Statistics\\
UC Berkeley\\
Berkeley, CA 94720, USA.}
\email{jding@stat.berkeley.edu}
\urladdr{}

\address{Eyal Lubetzky\hfill\break
Microsoft Research\\
One Microsoft Way\\
Redmond, WA 98052-6399, USA.}
\email{eyal@microsoft.com}
\urladdr{}

\address{Yuval Peres\hfill\break
Microsoft Research\\
One Microsoft Way\\
Redmond, WA 98052-6399, USA.}
\email{peres@microsoft.com}
\urladdr{}

\begin{abstract}
In the heat-bath Glauber dynamics for the Ising model
on the lattice, physicists believe that the spectral gap of the continuous-time
chain exhibits the following behavior.
For some critical inverse-temperature $\beta_c$, the inverse-gap is $O(1)$ for
$\beta < \beta_c$, polynomial in the surface area for $\beta = \beta_c$ and exponential in it
for $\beta > \beta_c$. This has been proved for $\Z^2$ except at criticality. So far, the only
underlying geometry where the critical behavior has been confirmed is the complete graph.
Recently, the dynamics for the Ising model on a regular tree, also known as the Bethe lattice,
has been intensively studied. The facts that the inverse-gap is bounded for $\beta < \beta_c$ and exponential
for $\beta > \beta_c$ were established, where $\beta_c$ is the critical spin-glass parameter,
and the tree-height $h$ plays the role of the surface area.

In this work, we complete the picture for the inverse-gap of the Ising model on the $b$-ary tree,
by showing that it is indeed polynomial in $h$ at criticality. The degree of our polynomial bound does not depend on $b$,
and furthermore, this result holds under any boundary condition. We also obtain
analogous bounds for the mixing-time of the chain.
In addition, we study the near critical behavior, and show that for $\beta > \beta_c$,
the inverse-gap and mixing-time are both $\exp[\Theta((\beta-\beta_c) h)]$.
\end{abstract}
\maketitle

\vspace{-0.175in}

\section{Introduction}

The \emph{Ising Model} on a finite graph $G=(V,E)$ with inverse-temperature $\beta \geq 0$ and no external magnetic field is defined as follows. Its set of possible \emph{configurations} is $\Omega = \{\pm1\}^V$, where each $\sigma\in\Omega$ assigns positive or negative \emph{spins} to the vertices of $G$. In the free boundary case, the probability that the system is at a given configuration $\sigma$ is given by the \emph{Gibbs distribution} $$\mu_G(\sigma)  = \frac1{Z(\beta)} \exp\Big(\beta \sum_{xy\in E}\sigma(x)\sigma(y)\Big)~,$$
where $Z(\beta)$ is the partition function. In the presence of a boundary condition $\tau\in \{\pm1\}^{\partial V}$ (that fixes the spins of some subset $\partial V\subset V$ of the sites), we let $\mu_G^\tau(\sigma)$ denote the Gibbs measure conditioned on $\sigma$ agreeing with $\tau$ on $\partial V$.

The heat-bath Glauber dynamics for the Ising model on $G$ is the Markov chain with the following transition rule: At each step,
a vertex is chosen uniformly at random, and its spin is updated according to $\mu^\tau_G$ conditioned on the spins of all the other vertices.
It is easy to verify that this chain is reversible with respect to the Gibbs distribution $\mu^\tau_G$. The continuous-time version of the dynamics associates each site with an independent Poisson clock of unit rate, determining the update times of this site as above (note that the continuous dynamics is $|V|$ times faster than the discrete dynamics).

The \emph{spectral-gap} of a reversible discrete-time chain, denoted by $\gap$, is $1-\lambda$, where $\lambda$ is the largest nontrivial eigenvalue of its transition kernel. The spectral-gap of the continuous-time process is defined analogously via the spectrum of its generator, and in the special case of Glauber dynamics for $\mu_G^\tau$, this gap is precisely $|V|$ times the discrete-time gap. This parameter
governs the rate of convergence to equilibrium in $L^2(\mu_G^\tau)$.

In the classical Ising model, the underlying geometry is the $d$-dimensional lattice, and there is a critical inverse-temperature $\beta_c$ where
the static Gibbs measure exhibits a phase transition with respect to long-range correlations between spins. While the main focus of
the physics community is on critical behavior (see the 20 volumes of \cite{DL}),
so far, most of the rigorous mathematical analysis was confined to
the non-critical regimes.

Supported by many experiments and studies in the theory of dynamical critical phenomena, physicists believe that the spectral-gap of the continuous-time dynamics on lattices has the following \emph{critical slowing down} behavior (e.g., \cite{HH,LF,Martinelli97,WH}):
At high temperatures ($\beta < \beta_c$) the inverse-gap is $O(1)$, at the critical $\beta_c$ it is polynomial in the surface-area and at low temperatures it is exponential in it. This is known for $\Z^2$ except at the critical $\beta_c$, and establishing the order of the gap at criticality seems extremely challenging. In fact, the only underlying geometry, where the critical behavior of the spectral-gap has been fully established, is the complete graph (see \cite{DLP}).

The important case of the Ising model on a regular tree, known as the Bethe lattice, has been intensively studied (e.g., \cite{BKMP,BRZ,CCCST1,CCCST2,CCST,EKPS,Ioffe1,Ioffe2,Lyons,MSW,PP}). On this canonical example of a \emph{non-amenable} graph (one whose boundary is proportional to its volume), the model exhibits a rich behavior. For example, it has two distinct critical inverse-temperatures: one for  uniqueness of the Gibbs state, and another for the purity of the free-boundary state. The latter, $\beta_c$, coincides with the critical spin-glass parameter.

As we later describe, previous results on the Ising model on a regular tree imply that the correct parameter to play the role of the surface-area is the tree-height $h$: It was shown that the inverse-gap is $O(1)$ for $\beta < \beta_c$ and exponential in $h$ for $\beta > \beta_c$, yet its critical behavior remained unknown.

In this work, we complete the picture of the spectral-gap of the dynamics for the critical Ising model on a regular tree, by establishing that it is
indeed polynomial in $h$. Furthermore, this holds under any boundary condition, and an analogous result is obtained for the $L^1$ (total-variation) mixing time, denoted by $\tmix$ (formally defined in Subsection \ref{subsec:prelim-tv}).

\begin{maintheorem}\label{thm-1}
Fix $b\geq 2$ and let $\beta_c = \arctanh(1/\sqrt{b})$ denote the critical inverse-temperature for the Ising model on the $b$-ary tree of height $h$.
Then there exists some constant $c>0$ independent of $b$, so that the following holds:
For any boundary condition $\tau$, the continuous-time Glauber dynamics for the above critical Ising model satisfies
$\gap^{-1} \leq \tmix = O(h^c)$.
\end{maintheorem}

One of the main obstacles in proving the above result is the arbitrary boundary condition, due to which the spin system loses its symmetry (and the task of analyzing the dynamics becomes considerably more involved). Note that, although boundary conditions are believed to only accelerate the mixing of the dynamics, even tracking
the effect of the (symmetric) all-plus boundary on lattices for $\beta>\beta_c$ is a formidable
open problem (see \cite{Martinelli94}).

In light of the above theorem and the known fact that the inverse-gap is exponential in $h$ at low temperatures ($\beta > \beta_c$ fixed), it
is natural to ask how the transition between these two phases occurs, and in particular, what the critical exponent of $\beta-\beta_c$ is. This is answered by the following theorem, which establishes that $\log(\gap^{-1})\asymp (\beta-\beta_c)h + \log h$ for small $\beta - \beta_c$. Moreover, this result also holds for $\beta = \beta_c + o(1)$, and thus pinpoints the transition to a polynomial inverse-gap at $\beta - \beta_c \asymp  \frac{\log h}{h}$.
\begin{maintheorem}\label{thm-2}
For some $\epsilon_0 > 0$, any $b \geq 2$ fixed and all $\beta_c < \beta < \beta_c + \epsilon_0$, where $\beta_c = \arctanh(1/\sqrt{b})$ is the critical spin-glass parameter, the following holds:
The continuous-time Glauber dynamics for the Ising model on a $b$-ary tree with inverse-temperature $\beta$ and free boundary satisfies
 \begin{equation}
  \label{eq-inverse-gap-transition}
\begin{array}{ll}
\gap^{-1} = h^{\Theta(1)} & \mbox{if }\beta = \beta_c + O(\frac{\log h}{h})~,\\
\gap^{-1} = \exp\left[\Theta\left((\beta-\beta_c)h\right)\right] & \mbox{otherwise.}
\end{array}
\end{equation}
Furthermore, both upper bounds hold under any boundary condition $\tau$, and \eqref{eq-inverse-gap-transition}
remains valid if $\gap^{-1}$ is replaced by $\tmix$.
\end{maintheorem}
In the above theorem and in what follows, the notation $f=\Theta(g)$ stands for $f=O(g)$ and $g=O(f)$.

Finally, our results include new lower bounds on the critical inverse-gap and the total-variation mixing-time (see Theorem~\ref{thm-tmix-trel-lower-bound}). The lower bound on $\gap^{-1}$ refutes a conjecture of \cite{BKMP}, according to which the continuous-time inverse-gap is linear in $h$. Our lower bound on $\tmix$ is of independent interest: Although in our setting the ratio between $\tmix$ and $\gap^{-1}$ is at most poly-logarithmic in $n$, the number of sites, we were able to provide a lower bound of order $\log n$ on this ratio without resorting to eigenfunction analysis.

\subsection{Background}\label{subsec:background}
The thoroughly studied question of whether the free boundary state is \emph{pure} (or \emph{extremal}) in the Ising model on the Bethe lattice can be formulated as follows: Does the effect that a typical boundary has
on the spin at the root vanish as the size of the tree tends to infinity?
It is well-known that one can sample a configuration for the tree according to the Gibbs distribution with free boundary by propagating spins along the tree (from a site to its children) with an appropriate bias (see Subsection \ref{subsec:prelim-ising} for details). Hence, the above question is equivalent to
asking wether the spin at the root can be reconstructed from its leaves, and as such has applications in Information Theory and Phylogeny (see \cite{EKPS} for further details).

In sufficiently high temperatures, there is a unique Gibbs state for the Ising model on a $b$-ary tree ($b \geq 2$), hence in particular the free boundary state is pure. The phase-transition with respect to the uniqueness of the Gibbs distribution occurs at the inverse-temperature
$\beta_u = \arctanh(1/b)$, as established in 1974 by Preston \cite{Preston}.

In \cite{CCCST2}, the authors studied the critical spin-glass model on the Bethe lattice (see also \cite{CCST,CCCST1}), i.e., the Ising model with a boundary of i.i.d.\ uniform spins. Following that work, it was finally shown in \cite{BRZ} that the phase-transition in the free-boundary extremality has the same critical inverse-temperature as in the spin-glass model,
$\beta_c = \arctanh(1/\sqrt{b})$. That is, the free-boundary state is pure iff $\beta \leq \beta_c$. This was later reproved in \cite{Ioffe1,Ioffe2}.

The inverse-gap of the Glauber dynamics for the Ising model on a graph $G$ was related in \cite{BKMP} to the \emph{cut-width} of the graph, $\xi(G)$, defined as follows: It is the minimum integer $m$, such that for some labeling of the vertices $\{v_1,\ldots,v_n\}$ and any $k\in[n]$, there are at most $m$ edges between $\{v_1,\ldots,v_k\}$ and $\{v_{k+1},\ldots,v_n\}$. The authors of \cite{BKMP} proved that for any bounded degree graph $G$, the continuous-time gap satisfies $\gap^{-1} = \exp[O(\xi(G)\beta)]$.

Recalling the aforementioned picture of the phase-transition of the gap, this supports the claim that the cut-width is the correct extension of the surface-area to general graphs.
One can easily verify that for $\Z_L^d$ (the $d$-dimensional box of side-length $L$) the cut-width has the same order as the surface-area $L^{d-1}$, while for a regular tree of height $h$ it is of order $h$.

Indeed, for the Ising model on a $b$-ary tree with $h$ levels and free boundary, it was shown in \cite{BKMP} that
the inverse-gap is $O(1)$ for all $\beta < \beta_c$, whereas for $\beta > \beta_c$ it satisfies $\log \gap^{-1} \asymp h$ (with constants that depend on $b$ and $\beta$). The behavior of the gap at criticality was left as an open problem: it is proved in \cite{BKMP} that
the critical $\gap^{-1}$ is at least linear in $h$ and conjectured that this is tight. A weaker conjecture of \cite{BKMP} states that $\gap^{-1}=\exp(o(h))$.

Further results on the dynamics were obtained in \cite{MSW}, showing that the log-Sobolev constant $\sob$ (defined in Section \ref{sec:prelim}) is uniformly bounded away from zero for $\beta < \beta_c$ in the free-boundary case, as well as for any $\beta$ under the all-plus boundary condition. While this implies that $\gap^{-1}=O(1)$ in these regimes, it sheds no new light on the behavior of the parameters $\gap,\sob$ in our setting of the critical Ising model on a regular tree with free-boundary.

\subsection{The critical inverse-gap and mixing-time}
Theorems~\ref{thm-1},\ref{thm-2}, stated above, establish that
on a regular tree of height $h$, the critical and near-critical continuous-time $\gap^{-1}$ and $\tmix$ are polynomial in $h$. In particular, this confirms the conjecture of \cite{BKMP} that the critical inverse-gap is $\exp(o(h))$.

Moreover, our upper bounds hold for any boundary condition, while matching the behavior of the free-boundary case:
Indeed, in this case the critical inverse-gap is polynomial in $h$ (as \cite{BKMP} showed it is at least linear),
and for $\beta-\beta_c > 0$ small we do have that $\log(\gap^{-1}) \asymp (\beta-\beta_c)h$.
For comparison, recall that under the all-plus boundary condition, \cite{MSW} showed that $\gap^{-1}=O(1)$ at all temperatures.

We next address the conjecture of \cite{BKMP} that the critical inverse-gap is in fact linear in $h$.
The proof that the critical $\gap^{-1}$ has order at least $h$ uses the same argument
that gives a tight lower bound at high temperatures: Applying
the Dirichlet form (see Subsection \ref{subsec:prelim-gap-sob})
to the sum of spins at the leaves as a test-function. Hence, the idea behind the above conjecture
is that the sum of spins at the boundary (that can be thought of as the magnetization) approximates the second eigenfunction also for $\beta = \beta_c$.

The following theorem refutes this conjecture, and en-route also implies that $\sob=o(1)$ at criticality.
In addition, this theorem provides a nontrivial lower bound on $\tmix$ that separates it from $\gap^{-1}$
(thus far, our bounds in Theorems~\ref{thm-1},\ref{thm-2} applied to both parameters as one).

\begin{maintheorem}\label{thm-tmix-trel-lower-bound}
Fix $b\geq 2$ and let $\beta_c = \arctanh(1/\sqrt{b})$ be the critical inverse-temperature for the Ising model on the $b$-ary tree with $n$ vertices. Then the corresponding discrete-time Glauber dynamics with free boundary satisfies:
\begin{align}
\gap^{-1} &\geq c_1 \, n \left(\log n\right)^{2}~,\label{eq-trel-lower-bound}\\
\tmix &\geq c_2 \, n \left(\log n\right)^{3}~,\label{eq-tmix-lower-bound}
\end{align}
for some $c_1,c_2 > 0$. Furthermore, $\tmix \geq c\, \gap^{-1} \log n$ for some $c > 0$.
\end{maintheorem}

Indeed, the above theorem implies that in continuous-time, $\gap^{-1}$ has order at least $h^2$
and $\tmix$ has order at-least $h^3$, where $h$ is again the height of the tree. By known facts on the log-Sobolev
constant (see Section \ref{sec:prelim}, Corollary~\ref{cor-tmix-trel-log-2}), in our setting we have $\tmix = O(\sob^{-1} h)$,
and it then follows that $\sob = O(h^{-2}) = o(1)$.

We note that by related results on the log-Sobolev constant, it follows that in the Ising model on a regular tree, for
any temperature and with any boundary condition we have
$\tmix = O( \gap^{-1} \log^2 n)$. In light of this, establishing a lower bound of order $\log n$ on the ratio between $\tmix$ and $\gap^{-1}$
is quite delicate (e.g., proving such a bound usually involves constructing a distinguishing statistic via a suitable eigenfunction (Wilson's method \cite{Wilson})).

\subsection{Techniques and proof ideas}

To prove the main theorem, our general approach is a recursive analysis of the spectral-gap via an appropriate block-dynamics (roughly put, multiple sites comprising a block
are updated simultaneously in each step of this dynamics; see Subsection~\ref{subsec:prelim-blocks} for a formal definition). This provides an estimate of the spectral-gap of the single-site dynamics in terms of those of the individual blocks and the block-dynamics chain itself (see \cite{Martinelli97}).
However, as opposed to most applications of the block-dynamics method,
where the blocks are of relatively small size, in our setting we must partition a tree of height $h$ to subtrees of height linear in $h$.
This imposes arbitrary boundary conditions on the individual blocks, and highly complicates the analysis of the block-dynamics chain.

In order to estimate the gap of the block-dynamics chain, we apply the method of Decomposition of Markov chains, introduced in \cite{JSTV}
(for details on this method see Subsection~\ref{subsec:prelim-decomposition-mc}). Combining this method with a few other ideas
(such as establishing contraction and controlling the external field in certain chains), the proof of Theorem~\ref{thm-1} is reduced into the following
spatial-mixing/reconstruction type problem. Consider the procedure, where we assign the spins of the boundary given the value at the root of the tree,
then reconstruct the root from the values at the boundary. The key quantity required in our setting is the difference in the expected
outcome of the root, comparing the cases where its initial spin was either positive or negative.

This quantity was studied by \cite{PP} in the free-boundary case, where it was related to capacity-type parameters of the tree (see \cite{EKPS} for a related result corresponding to the high temperature regime). Unfortunately, in our case we have an arbitrary boundary condition, imposed by the block-dynamics. This eliminates the symmetry of the system, which was a crucial part of the arguments of \cite{PP}. The most delicate step in the proof of Theorem~\ref{thm-1} is the extension of these results of \cite{PP} to any boundary condition. This is achieved by carefully tracking down the effect of the boundary on the expected reconstruction result in each site, combined with correlation inequalities and an analytical study of the corresponding log-likelihood-ratio function.

The lower bound on the critical inverse-gap reflects the change in the structure of the dominant eigenfunctions between high and low temperatures.
At high temperatures, the sum of spins on the boundary gives the correct order of the gap. At low temperatures, a useful lower bound on $\gap^{-1}$ was shown in \cite{BKMP} via the recursive-majority function (intuitively, this reflects the behavior at the root: Although this spin may occasionally flip its value, at low temperature it constantly tends to revert to its biased state).
Our results show that at criticality, a lower bound improving upon that of \cite{BKMP} is obtained by essentially merging the above two functions into
a weighted sum of spins, where the weight of a spin is determined by its tree level.

To establish a lower bound on $\tmix$ of order $\gap^{-1}h$, we consider a certain \emph{speed-up} version of the dynamics: a block-dynamics,
whose blocks are a mixture of singletons and large subtrees. The key ingredient here is the Censoring Inequality of Peres and Winkler \cite{PW},
that shows that this dynamics indeed mixes as fast as the usual (single-site) one. We then consider a series of modified versions
of this dynamics, and study their mixing with respect to the total-variation and Hellinger distances. In the end, we arrive at a product chain,
whose components are each the single-site dynamics on a subtree of height linear in $h$. This latter chain provides the required lower bound on $\tmix$.

%

\subsection{Organization}
The rest of this paper is organized as follows. Section~\ref{sec:prelim} contains several preliminary facts and definitions.
In Section~\ref{sec:spatial-mixing} we prove a spatial-mixing type result on the critical and near-critical Ising model on a tree with an arbitrary boundary condition. This then serves as one of the key ingredients in the proof of the main result, Theorem~\ref{thm-1}, which appears in Section~\ref{sec:thm-1}. In Section~\ref{sec:lower-bound} we prove Theorem~\ref{thm-tmix-trel-lower-bound}, providing
the lower bounds for the critical inverse-gap and mixing-time. Section~\ref{sec:phase-transition} contains the proof of Theorem~\ref{thm-2}, addressing
the near-critical behavior of $\gap^{-1}$ and $\tmix$. The final section, Section~\ref{sec:conclusion},
 is devoted to concluding remarks and open problems.

\section{Preliminaries}\label{sec:prelim}
\subsection{Total-variation mixing}\label{subsec:prelim-tv}
Let $(X_t)$ be an aperiodic irreducible Markov chain on a finite state space $\Omega$, with stationary distribution $\pi$.
For any two distributions $\phi,\psi$ on $\Omega$, the \emph{total-variation distance} of $\phi$ and $\psi$ is defined as
$$\|\phi-\psi\|_\mathrm{TV} \deq \sup_{A \subset\Omega} \left|\phi(A) - \psi(A)\right| = \frac{1}{2}\sum_{x\in\Omega} |\phi(x)-\psi(x)|~.$$
 The (worst-case) total-variation \emph{mixing-time} of $(X_t)$, denoted by $\tmix(\epsilon)$ for $0 < \epsilon < 1$, is defined to be $$ \tmix(\epsilon) \deq \min\Big\{t : \max_{x \in \Omega} \| \P_x(X_t \in \cdot)- \pi\|_\mathrm{TV} \leq \epsilon \Big\}~,$$
where $\P_x$ denotes the probability given that $X_0=x$. As it is easy and well known (cf., e.g., \cite{AF}*{Chapter 4}) that
the spectral-gap of $(X_t)$ satisfies $ \gap^{-1} \leq \tmix\left(1/\mathrm{e}\right)$, it will be convenient to use the abbreviation
$$\tmix \deq \tmix\left(1/\mathrm{e}\right)~.$$
Analogously, for a continuous-time chain on $\Omega$ with heat-kernel $H_t$, we define $\tmix$ as
the minimum $t$ such that $\max_{x \in \Omega} \| H_t(x,\cdot)- \pi\|_\mathrm{TV} \leq 1/\mathrm{e}$.

\subsection{The Ising model on trees}\label{subsec:prelim-ising}
When the underlying geometry of the Ising model is a tree with free boundary condition, the Gibbs measure has a natural constructive representation.
This appears in the following well known claim (see, e.g., \cite{EKPS} for more details).
\begin{claim}\label{claim-Ising-representation}
Consider the Ising model on a tree $T$ rooted at $\rho$ with free boundary condition and at the
inverse-temperature $\beta$. For all $e\in E(T)$, let $\eta_e \in \{\pm1\}$ be i.i.d. random variables with $\P(\eta_e = 1)
= (1+\tanh \beta)/2$. Furthermore, let $\sigma(\rho)$ be a uniform spin, independent of $\{\eta_e\}$, and for $v \neq \rho$,
$$\sigma(v) = \sigma(\rho) \prod_{e \in \mathcal{P}(\rho, v)} \eta_e~,\quad \mbox{where $\mathcal{P}(\rho, v)$ is the simple path from $\rho$ to $v$.}$$
Then the distribution of the resulting $\sigma$ is the corresponding Gibbs measure.
\end{claim}
In light of the above claim, one is able to sample a configuration according to Gibbs distribution on a tree with free boundary condition using the following
simple scheme: Assign a uniform spin at the root $\rho$, then scan the tree from top to bottom, successively assigning each site with a spin according to the value at its parent. More precisely, a vertex is assigned the same spin as its parent with probability $(1 + \tanh \beta)/2$, and the opposite one otherwise. Equivalently, a vertex inherits the spin of its parent with probability $\tanh \beta$, and otherwise it receives an independent uniform spin. Finally, for the conditional Gibbs distribution given a plus spin at the root $\rho$, we assign $\rho$ a plus spin rather than a uniform spin, and carry on as above.

However, notice that the above \emph{does not} hold for the Ising model in the presence of a boundary condition,
which may impose a different external influence on different sites.

\subsection{$L^2$-capacity}\label{subsec:prelim-capacity}
The authors of \cite{PP} studied certain spatial mixing properties of the Ising model on trees (with free or all-plus boundary conditions), and related them to
the $L^p$-capacity of the underlying tree. In Section \ref{sec:spatial-mixing}, we extend some of the results of \cite{PP}
to the (highly asymmetric) case of a tree with an arbitrary boundary condition, and relate a certain ``decay of correlation''
 property to the $L^2$-capacity of the tree, defined as follows.

Let $T$ be a tree rooted at $\rho$, denote its leaves by $\partial T$, and throughout the paper, write $(u,v)\in E(T)$ for the directed edge between a vertex $u$ and its child $v$. We further define $\dist(u,v)$ as the length (in edges) of the simple path connecting $u$ and $v$ in $T$.

For each $e\in E(T)$, assign the \emph{resistance} $R_e \geq 0 $ to
the edge $e$. We say that a non-negative function $f : E(T)\to\R$ is a \emph{flow} on $T$ if the following holds
for all $(u, v) \in E(T)$ with $v \not\in \partial T$:
$$f(u, v) = \sum_{(v, w) \in E(T)} f(v, w)~,$$
that is, the incoming flow equals the outgoing flow on each internal vertex $v$ in $T$.
For any flow $f$, define its \emph{strength} $|f|$ and \emph{voltage} $V(f)$ by
\begin{align*}
|f| \deq \sum_{(\rho, v)\in E(T)} f(\rho, v)~,~V(f) \deq \sup\bigg\{\sum_{e\in \mathcal{P}(\rho, w)} f(e)R_e: w \in \partial T\bigg\}~,
\end{align*}
where $\mathcal{P}(\rho, w)$ denotes the simple path from $\rho$ to $w$.
Given these definitions, we now define the $L^2$-capacity $\mathrm{cap}_2(T)$ to be
$$\mathrm{cap}_2(T) \deq \sup\{|f|: f \mbox{ is a flow with } V(f) \leq 1\}~.$$
For results on the $L^2$-capacity of general networks (and more generally, $L^p$-capacities, where the expression $f(e)R_e$ in the above definition of $V(f)$ is replaced by its $(p-1)$ power), as part of Discrete Nonlinear Potential Theory,
cf., e.g., \cite{MY}, \cite{Soardi93}, \cite{Soardi94} and the references therein.

For our proofs, we will use the well-known fact that the $L^2$-capacity of the tree $T$ is precisely the \emph{effective conductance} between the root $\rho$ and the leaves $\partial T$, denoted by $C_\mathrm{eff}(\rho \leftrightarrow \partial T)$. See, e.g., \cite{LP} for further information on electrical networks.

\subsection{Spectral gap and log-Sobolev constant}\label{subsec:prelim-gap-sob}
Our bound on the mixing time of Glauber dynamics for the Ising model on trees will be derived from a recursive analysis of the spectral gap of this chain. This analysis uses spatial-mixing type results (and their relation to the above mentioned $L^2$~capacity) as a building block. We next describe how the mixing-time can be bounded via the spectral-gap in our setting.

The spectral gap and log-Sobolev constant of a reversible Markov chain with stationary distribution $\pi$ are given by the following Dirichlet forms (see, e.g., \cite{AF}*{Chapter 3,8}):
\begin{align}\label{eq-dirichlet-form}
\gap = \inf_f \frac{\mathcal{E}(f)}
{\var_{\pi}(f)}~,~\quad \sob = \inf_{f} \frac{\mathcal{E}(f)}{\operatorname{Ent}(f)}~,
\end{align}
where
\begin{align}\label{eq-dirichlet-form-e}
&\mathcal{E}(f) = \left<(I-P)f,f\right>_{\pi} = \frac{1}{2}\sum_{x,y\in\Omega}\left[f(x)-f(y)\right]^2\pi(x)P(x,y)~,\\
&\operatorname{Ent}_\pi(f) = \E_\pi\left(f^2 \log ( f^2 /\E_\pi f^2)\right) ~.
\end{align}
By bounding the log-Sobolev constant, one may obtain remarkably sharp upper bounds on the $L^2$ mixing-time: cf., e.g., \cite{DS1,DS2,DS,DS3,SaloffCoste}. The following result
of Diaconis and Saloff-Coste \cite{DS}*{Theorem 3.7} (see also \cite{SaloffCoste}*{Corollary~2.2.7}) demonstrates this powerful
method; its next formulation for discrete-time appears in \cite{AF}*{Chapter 8}. As we are interested in total-variation mixing,
we write this bound in terms of $\tmix$, though it in fact holds also for the (larger) $L^2$ mixing-time.
\begin{theorem}[\cite{DS}, \cite{SaloffCoste}, reformulated]\label{thm-logsolev-mixing}
For any reversible finite Markov chain with stationary distribution $\pi$,
$$\tmix(1/\mathrm{e}) \leq \frac14 \sob^{-1} \left( \log\log(1/\pi^*) + 4\right)~, $$
where $\pi^*=\min_x \pi(x)$.
\end{theorem}
We can then apply a result of \cite{MSW}, which provides a useful bound on $\sob$ in terms of $\gap$ in our setting, and obtain an upper bound on the mixing-time.
\begin{theorem}[\cite{MSW}*{Theorem 5.7}]\label{thm-gap-logsobolev-tree}
There exists some $c > 0$ such that the Ising model on the $b$-ary tree with $n$ vertices satisfies $\sob \geq c \cdot \gap /\log n$.
\end{theorem}
Note that the proof of the last theorem holds for any $\beta$ and under any boundary condition.
Combining Theorems \ref{thm-logsolev-mixing} and \ref{thm-gap-logsobolev-tree}, and noticing that $\pi_* \geq 2^{-n} \exp(-2\beta n)$
(as there are $2^n$ configurations, and the ratio between the maximum and minimum probability of a configuration is at most $\exp(2\beta n)$), we obtain the following corollary:
\begin{corollary}\label{cor-tmix-trel-log-2}
The Glauber dynamics for the Ising model on a $b$-ary tree with $n$ vertices satisfies
$\tmix = O\left( \sob^{-1} \log n\right) = O\left( \gap^{-1} \log^2 n\right)$ for any $\beta$ and any boundary condition.
\end{corollary}

The above corollary reduces the task of obtaining an upper bound for the mixing-time into establishing a suitable lower bound on the spectral gap. This will be achieved using a block dynamics analysis.

\subsection{From single site dynamics to block dynamics}\label{subsec:prelim-blocks}
Consider a cover of $V$ by a collection of subsets $\{B_1,\ldots,B_k\}$, which we will refer to as ``blocks''.
The \emph{block dynamics} corresponding to $B_1,\ldots,B_k$ is the Markov chain, where at each step a uniformly
chosen block is updated according to the stationary distribution given the rest of the system. That is,
the entire set of spins of the chosen block is updated simultaneously, whereas all other spins remain unchanged.
One can verify that the block dynamics is reversible with respect to the Gibbs distribution $\mu_n$.

Recall that, given a subset of the sites $U \subset V$, a boundary condition $\eta$ imposed on $U$ is
the restriction of the sites $U^c = V \setminus U$ to all agree with $\eta$ throughout the dynamics, i.e., only sites in $U$ are considered for updates. It will sometimes be useful to consider $\eta \in \Omega$ (rather than a configuration of the sites $U^c$),
in which case only its restriction to $U^c$ is accounted for.


The following theorem shows the useful connection between the single-site dynamics and the block dynamics.
This theorem appears in \cite{Martinelli97} in a more general setting, and following is its reformulation for
the special case of Glauber dynamics for the Ising model on a finite graph
with an arbitrary boundary condition. Though the original theorem is stated for the continuous-time
dynamics, its proof naturally extends to the discrete-time case; we provide its details for completeness.
\begin{proposition}[\cite{Martinelli97}*{Proposition 3.4}, restated]\label{prop-block-single}
Consider the discrete time Glauber dynamics on a $b$-ary tree with boundary condition $\eta$.
Let $\gap_U^\eta$ be the spectral-gap of the single-site
dynamics on a subset $U\subset V$ of the sites, and $\gap_\mathcal{B}^\eta$ be the spectral-gap of the block dynamics corresponding to $B_1,\ldots,B_k$, an arbitrary cover of a vertex set $W\subset V$.  The following holds:
\begin{equation*}
  \gap_W^\eta \geq \frac{k}{|W|} \gap_\mathcal{B}^\eta \inf_i \inf_{\varphi} |B_i| \gap_{B_i}^\varphi \left(\sup_{x\in W} \#\{i:B_i \ni x\}\right)^{-1}~.
\end{equation*}
\end{proposition}
\begin{proof} Let $P$ denote the transition kernel of the above Glauber dynamics. Defining
$$g \deq \inf_i \inf_\varphi |B_i| \gap^\varphi_{B_i}~,$$
the Dirichlet form \eqref{eq-dirichlet-form} gives that, for any function $f$,
$$\var^\varphi_{B_i}(f) \leq \frac{\mathcal{E}_{B_i}^\varphi(f)}{\gap^\varphi_{B_i}}  \leq \frac{|B_i|}{g} \mathcal{E}_{B_i}^\varphi(f)~.$$
Combining this with definition \eqref{eq-dirichlet-form-e} of $\mathcal{E}(\cdot)$,
\begin{align*}
\mathcal{E}^\eta_\mathcal{B}(f) = \frac{1}{k}\sum_{\varphi \in \Omega} \mu_W^\eta(\varphi) \sum_i \var^\varphi_{B_i}(f) \leq \frac{1}{kg}\sum_{\varphi \in \Omega} \mu_W^\eta(\varphi) \sum_i |B_i| \mathcal{E}_{B_i}^\varphi(f)~.
\end{align*}
On the other hand, definition \eqref{eq-dirichlet-form-e} again implies that
\begin{align*}
&\sum_{\varphi \in \Omega} \mu_W^\eta(\varphi)  \sum_{i} |B_i| \mathcal{E}_{B_i}^\varphi(f) \\
&= \sum_{\varphi \in \Omega} \mu_W^\eta(\varphi)\frac{1}{2} \sum_{\sigma\in \Omega} \sum_i \sum_{x\in B_i} \mu^\varphi_{B_i}(\sigma)|B_i|P_{B_i}^\varphi(\sigma, \sigma^x)[f(\sigma) - f(\sigma^x)]^2\\
 &\leq \frac{1}{2}\sup_{x\in W} \#\{i:B_i \ni x\} \sum_{\sigma \in \Omega} \mu_W^\eta(\sigma)\sum_{x\in W} |W|P^\sigma_{W}(\sigma, \sigma^x)[f(\sigma) - f(\sigma^x)]^2\\
 &=|W|\sup_{x\in W}\#\{i:B_i \ni x\} \mathcal{E}^\eta_W(f) ~,
\end{align*}
where $\sigma^x$ is the configuration obtained from $\sigma$ by flipping the spin at $x$, and we used the fact that
$$|B_i|P^\sigma_{B_i}(\sigma, \sigma^x) = |W| P_W^\sigma(\sigma, \sigma^x)\quad\mbox{ for any $i\in [k]$ and $x\in B_i$}~.$$
Altogether, we obtain that
$$\mathcal{E}^\eta_\mathcal{B}(f) \leq \frac{|W|}{kg}\sup_{x\in W}\#\{i:B_i \ni x\} \mathcal{E}^\eta_W(f)~.$$
Recalling that the single-site dynamics and the block-dynamics have the same stationary measure,
$$ \frac{\mathcal{E}^\eta_\mathcal{B}(f)}{\var_W^\eta(f)}=
\frac{\mathcal{E}^\eta_\mathcal{B}(f)}{\var_\mathcal{B}^\eta(f)} \geq \gap_\mathcal{B}^\eta$$
(where we again applied inequality \eqref{eq-dirichlet-form}), thus
$$\frac{\mathcal{E}^\eta_W(f)}{\var_W^\eta(f)} \geq \frac{k}{|W|}g \Big(\sup_{x\in W}\#\{i:B_i \ni x\} \Big)^{-1}\gap^\eta_\mathcal{B} ~.$$
The proof is now completed by choosing $f$ to be the eigenfunction that corresponds to the second eigenvalue of $P_W^\eta$ (achieving $\gap_W^\eta$), with a final application of \eqref{eq-dirichlet-form}.
\end{proof}

The above proposition can be applied, as part of the spectral gap analysis, to reduce the size of the base graph (though with an arbitrary boundary condition), provided that one can estimate the gap of the corresponding block dynamics chain.

\subsection{Decomposition of Markov chains}\label{subsec:prelim-decomposition-mc}
In order to bound the spectral gap of the block dynamics, we require a result of \cite{JSTV}, which analyzes the spectral gap of a Markov chain via its decomposition into a {\em projection} chain and a {\em restriction} chain.

Consider an ergodic Markov chain on a finite state space $\Omega$ with transition kernel $P : \Omega \times \Omega \to [0,1]$ and stationary distribution $\pi: \Omega \to [0,1]$. We assume that the Markov chain is \emph{time-reversible}, that is to say, it satisfies the detailed balance condition
$$\pi(x) P(x, y)  = \pi(y)P(y,x) \mbox{ for all $x, y \in \Omega$}~.$$
Let $\Omega = \Omega_0 \cup \ldots \cup \Omega_{m-1}$ be a decomposition of the state space into $m$ disjoint sets. Writing $[m] \deq \{0, \ldots, m-1\}$, we define $\bar{\pi} : [m] \to [0,1]$ as
$$\bar{\pi}(i) \deq \pi(\Omega_i) = \sum_{x \in \Omega_i} \pi (x)$$
and define $\bar{P}: [m] \times [m] \to [0,1]$ to be
$$\bar{P}(i,j) \deq \frac1{\bar{\pi}(i)} \sum_{x \in \Omega_i, y\in \Omega_j} \pi(x) P(x, y)~.$$
The Markov chain on the state space $[m]$ whose transition kernel is $\bar{P}$ is called the \emph{projection} chain, induced by the partition $\Omega_0,\ldots,\Omega_{m-1}$. It is easy to verify that, as the original Markov chain is reversible with respect to $\pi$, the projection chain is reversible with respect to the stationary distribution $\bar{\pi}$.

In addition, each $\Omega_i$ induces a \emph{restriction} chain, whose transition kernel $P_i : \Omega_i \times \Omega_i \to [0, 1]$ is given by
$$P_i(x, y) =\begin{cases}
P(x, y), &\mbox{if } x\neq y,\\
1 - \sum_{z \in \Omega_i \backslash \{x\}} P(x, z), &\mbox{if } x=y~.\end{cases}$$
Again, the restriction chain inherits its reversibility from the original chain,
and has a stationary measure $\pi_i$, which is simply $\pi$ restricted to $\Omega_i$: $$\pi_i(x) \deq\pi(x)/\bar{\pi}(i)~\mbox{ for all $x\in\Omega_i$}~.$$
In most applications, the projection chain and the different restriction
chains are all irreducible, and thus the various stationary distributions
$\bar{\pi}$ and $\pi_0, \ldots ,\pi_{m-1}$ are all unique.

The following result provides a lower bound on the spectral gap of the original Markov chain given its above described decomposition:
\begin{theorem}[\cite{JSTV}*{Theorem 1}]\label{thm-gap-decomposition}
Let $P$ be the transition kernel of a finite reversible Markov chain, and let $\gap$ denote its spectral gap.
Consider the decomposition of the chain into a projection chain and $m$ restriction chains, and denote their corresponding spectral gaps by $\bar{\gap}$ and $\gap_0,\ldots,\gap_{m-1}$. Define
 $$\gap_{\mathrm{min}} \deq \min_{i\in [m]} \gap_i~,~\quad\gamma \deq \max_{i\in [m]} \max_{x\in \Omega_i} \sum_{y\in \Omega\backslash \Omega_i} P(x, y)~.$$
Then $\gap$, the spectral gap of the original Markov chain, satisfies:
$$\gap \geq \frac{\bar{\gap}}{3}\;\wedge\; \frac{\bar{\gap}\cdot \gap_{\min}}{3\gamma + \bar{\gap}}~.$$
\end{theorem}
The main part of Section \ref{sec:thm-1} will be devoted to the analysis of the projection chain, in an effort to bound the spectral gap of our block dynamics via the above theorem.

\section{Spatial mixing of Ising model on trees}\label{sec:spatial-mixing}
In this section, we will establish a spatial-mixing type result for the Ising model on a general (not necessarily regular) finite tree under an arbitrary boundary condition. This result (namely, Proposition \ref{prop-delta}) will later serve as the main ingredient in the proof of Theorem \ref{thm-1} (see Section \ref{sec:thm-1}).
Throughout this section, let $\beta > 0$ be an arbitrary inverse-temperature and $\theta = \tanh \beta$.

We begin with a few notations. Let $T$ be a tree rooted at $\rho$ with a boundary condition $\tau\in\{\pm1\}^{\partial T}$ on its leaves, and
$\mu^\tau$ be the corresponding Gibbs measure.

For any $v \in T$, denote by $T_v$ the subtree of $T$ containing $v$ and its all descendants. In addition, for any
$B \subset A \subset T$ and $\sigma \in \{\pm1\}^A$, denote by $\sigma_B$ the restriction of $\sigma$ to the sites of $B$.
We then write $\mu_v^\tau$ for the Gibbs measure on the subtree $T_v$ given the boundary $\tau_{\partial T_v}$.

Consider $\hat{T} \subset T\setminus\partial T $, a subtree of $T$ that contains the root $\rho$, and write $\hat{T}_v = T_v \cap \hat{T}$.
Similar to the above definitions for $T$, we
denote by $\hat{\mu}^\xi$ the Gibbs measure on $\hat{T}$ given the boundary condition $\xi\in \{\pm1\}^{\partial \hat{T}}$, and
 let $\hat{\mu}_v^\xi$ be the Gibbs measure on $\hat{T}_v$ given the boundary $\xi_{\partial \hat{T}_v}$.

The following two measures are the conditional distributions of $\mu_v^\tau$ on the boundary of the subtree $\hat{T}_v$ given the spin at its
root $v$:
\begin{align*}
Q_v^+ (\xi )&\deq \mu_v^\tau \big( \sigma_{\partial\hat{T}_v} = \xi_{\partial\hat{T}_v} \mid \sigma(v) = 1\big)& \mbox{ for }\xi\in \{\pm1\}^{\partial \hat{T}}~,\\
Q_v^- (\xi )&\deq \mu_v^\tau \big( \sigma_{\partial\hat{T}_v} = \xi_{\partial\hat{T}_v} \mid \sigma(v) = -1\big)& \mbox{ for }\xi\in \{\pm1\}^{\partial \hat{T}}~.\\
\end{align*}
We can now state the main result of this section, which addresses the problem of reconstructing
the spin at the root of the tree from its boundary.
\begin{proposition}\label{prop-delta}
Let $\hat{T}$ be as above, let $0 < \theta \leq \frac34$ and define
\begin{equation*}\Delta \deq \int \hat{\mu}^{\xi}(\sigma(\rho) = 1) d Q_\rho^+(\xi) - \int \hat{\mu}^{\xi}(\sigma(\rho) = 1) d Q_\rho^-(\xi)~.
\end{equation*}
Then there exists an absolute constant $\kappa > \frac1{100}$ such that
 $$\Delta \leq \frac{\mathrm{cap}_2 (\hat{T})}{\kappa (1-\theta)}~,$$
where the resistances are assigned to be $R_{(u,v)} = \theta^{-2\dist(\rho,v)}$.
Furthermore, this also holds for any external field $h \in \R$ on the root $\rho$.
\end{proposition}
To prove the above theorem, we consider the notion of the log-likelihood ratio at a vertex $v$ with respect to $\hat{T}_v$ given the boundary $\xi_{\partial{\hat{T}}_v}$:
\begin{equation}\label{eq-xv-def}
x_v^\xi = \log\bigg(\frac{\hat{\mu}_v^\xi (\sigma(v) = +1)}{\hat{\mu}_v^\xi (\sigma(v) = -1)}\bigg)~,
\end{equation}
as well as the following quantity, analogous to $\Delta$ (defined in Proposition \ref{prop-delta}):
\begin{equation}
  \label{eq-mv-def}
 m_v \deq \int x_v^\xi d Q_v^+ - \int x_v^\xi d Q_v^-~.
\end{equation}
As we will later explain, $m_v \geq 0$ for any $v \in T$, and we seek an upper bound on this quantity.
One of the main results of \cite{PP} was such an estimate for the case of \emph{free} boundary condition,
yet in our setting we have an \emph{arbitrary} boundary condition (adding a considerable amount
of difficulties to the analysis). The following theorem
extends the upper bound on $m_\rho$ to any boundary; to avoid confusion, we formulate this bound in terms
 of the same absolute constant $\kappa$ given in Proposition \ref{prop-delta}.
\begin{theorem}\label{thm-log-likelihood-decay} Let $\hat{T}$ and $m_\rho$ be as above, and let $0 < \theta \leq \frac34$.
Then there exists an absolute constant $\kappa > \frac1{100}$ such that
$$m_\rho \leq \frac{\mathrm{cap}_2 (\hat{T})}{\kappa (1-\theta)/4}~,$$
where the resistances are assigned to be $R_{(u,v)} = \theta^{-2\dist(\rho,v)}$.
\end{theorem}

\subsection*{Proof of Theorem \ref{thm-log-likelihood-decay}}
As mentioned above, the novelty (and also the main challenge) in the result stated in Theorem \ref{thm-log-likelihood-decay} is the presence of the arbitrary boundary condition $\tau$, which eliminates most of the symmetry that one has in the free boundary case.
Note that this symmetry was a crucial ingredient in the proof of \cite{PP} for the free boundary case (namely, in that case $Q_v^+$ and $Q_v^-$ are naturally symmetric).

In order to tackle this obstacle, we need to track down the precise influence of the boundary condition $\tau$ on each vertex $v\in T$. We then incorporate this information in the recursive analysis that appeared (in a slightly different form) in \cite{Lyons}.
This enables us to relate the recursion relation of the $m_v$-s to that of the $L^2$-capacity.

The following quantity captures the above mentioned influence of $\tau$ on a given vertex $v\in T$:
\begin{equation}
  \label{eq-xv*-def}
x_v^* = \log \bigg(\frac{\mu_v^\tau(\sigma(v) = 1)}{\mu_v^\tau(\sigma(v) = -1)}\bigg)~.
\end{equation}
Notice that $x_v^*$ has a similar form to $x_v^\xi$ (defined in \eqref{eq-xv-def}), and
is essentially the log-likelihood ratio at $v$ induced by the boundary condition $\tau$. The quantity
$x_v^\xi$ is then the log-likelihood ratio that in addition considers the extra constraints imposed by $\xi$.
Also note that a free boundary condition corresponds to the case where $x_v^* = 0$ for all $v\in T$.

To witness the effect of $x_v^*$, consider the probabilities of propagating a spin from a parent $v$ to its child $w$,
formally defined by
\begin{align*}
p^\tau_{v, w}(1,1) &\deq \mu_{v}^\tau(\sigma(w) = 1 \mid \sigma(v) = 1)~,\\
 p^\tau_{v, w}(1,-1) &\deq \mu_{v}^\tau(\sigma(w) = -1 \mid\sigma(v) = 1)~;
\end{align*}
we define $p^\tau_{v, w}(-1,1)$ and $ p^\tau_{v, w}(-1,-1)$ analogously. The next simple lemma shows the relation between $x_v^*$ and these probabilities.

\begin{lemma}\label{lem-influence-tau}
The following holds for any $(v,w)\in T$:
$$p^\tau_{v, w}(1,1) - p^\tau_{v, w}(-1,1) = D_w^* \theta~,$$
where $D_w^* \deq (\cosh \beta)^2 / \left(\cosh^2 \beta + \cosh^2(x_w^*/2) -1\right)$.
\end{lemma}
\begin{proof}
Recalling definition \eqref{eq-xv*-def} of $x_v^*$, we can translate the boundary condition $\tau$ into an external field $x_w^*/2$ on the vertex $w$ when studying the distribution of its spin. Hence,
\begin{align*}
p^\tau_{v, w}(1,1) - p^\tau_{v, w}(-1,1) &= \frac{\mathrm{e}^{\beta + x_w^*/2}}{\mathrm{e}^{\beta + x_w^*/2} + \mathrm{e}^{-\beta - x_w^*/2}} -  \frac{\mathrm{e}^{-\beta + x_w^*/2}}{\mathrm{e}^{-\beta + x_w^*/2} + \mathrm{e}^{\beta - x_w^*/2}}\\
&=\frac{\mathrm{e}^{2\beta} - \mathrm{e}^{-2\beta}}{\mathrm{e}^{x_w^*} + \mathrm{e}^{-x_w^*} + \mathrm{e}^{2\beta} + \mathrm{e}^{-2\beta}}\\
&= \frac{\cosh^2 \beta}{\cosh^2 \beta + \cosh^2(x_w^*/2) -1} \tanh \beta~,
\end{align*}
as required.
\end{proof}
\begin{remark}
In the free boundary case, we have  $p_{v, w}(1,1) - p_{v, w}(-1,1) = \theta$.
For a boundary condition $\tau$, the coefficient $0<D_w^* \leq 1$ represents the contribution of this boundary to the propagation
probability.
\end{remark}

We now turn our attention to $m_v$. As mentioned before, the fact that $m_v \geq 0$ follows
from its definition \eqref{eq-mv-def}. Indeed, the monotonicity of the Ising model implies that the measure $Q_v^+$ stochastically dominates the measure $Q_v^-$ (with respect to the natural partial order on the configurations of $\partial\hat{T}_v$). For instance,
it is easy to see this by propagating $1$ and $-1$ spins from the root to the bottom, and applying a monotone coupling on these two processes. Finally, $x_v^\xi$ is monotone increasing in $\xi$ (again by the monotonicity of the Ising model), thus $m_v \geq 0$.

The first step in establishing the recursion relation of $m_v$ (that would lead to the desired upper bound)
would be to relate $m_v$ to some quantities associated with its children, as stated next.
\begin{lemma}\label{lem-recursion-m-f}
For any $v\in \hat{T}\setminus \partial\hat{T}$, we have that
$$m_v = \theta \sum_{w:(v,w)\in \hat{T}} D_w^*\Big(\int f(x_w^\xi) d Q_w^+(\xi) - \int f(X_w^\xi) d Q_w^-(\xi)\Big)~,$$
where
\begin{equation}\label{eq-def-f}f(x) = \log \bigg(\frac{\cosh(x/2) + \theta \sinh(x/2)}{\cosh(x/2) - \theta \sinh(x/2)}\bigg)~.\end{equation}
\end{lemma}
\begin{proof}
We need the following well-known lemma, that appeared in \cite{Lyons} in a slightly different form; see also \cite{Baxter} and \cite{PP}*{Lemma 4.1}.
\begin{lemma}[\cite{Lyons},\cite{PP} (reformulated)]
Let $f$ be as in \eqref{eq-def-f}. For all $v\in \hat{T}\setminus\partial\hat{T}$ and $\xi \in \{\pm1\}^{\partial\hat{T}}$, the following holds:
$x_v^\xi = \sum_{w:(v,w)\in \hat{T}} f(x_w^\xi)$.
\end{lemma}
According to this lemma, we obtain that
\begin{align}\label{eq-recursion-m_v}
m_v =\sum_{w:(v,w)\in\hat{T}}\Big( \int f(x_w^\xi) d Q_v^+(\xi) - \int f(x_w^\xi) d Q_v^-(\xi) \Big)~.
\end{align}
Noting that $x_w^\xi$ is actually a function of $\xi_{\partial\hat{T}_w}$, we get that
$$\int f(x_w^\xi) d Q_v^+(\xi) = \int f(x_w^\xi) d(p^\tau_{v, w}(1,1) Q_w^+(\xi) + p^\tau_{v, w}(1,-1) Q_w^-(\xi))~,$$
and similarly, we have
$$\int f(x_w^\xi) d Q_v^-(\xi) = \int f(x_w^\xi) d(p^\tau_{v, w}(-1,1) Q_w^+(\xi) + p^\tau_{v, w}(-1,-1) Q_w^-(\xi))~.$$
Combining these two equalities, we deduce that
\begin{align}\label{eq-bound-f-x}
&\int f(x_w^\xi) d Q_v^+(\xi)  - \int f(x_w^\xi) d Q_v^-(\xi) \nonumber\\
&= (p^\tau_{v, w}(1, 1) - p^\tau_{v,w}(-1,1))\Big(\int f(x_w^\xi) d Q_w^+(\xi) - \int f(x_w^\xi) d Q_w^-(\xi)\Big)\nonumber\\
&=  \theta D_w^*  \Big(\int f(x_w^\xi) d Q_w^+(\xi) - \int f(x_w^\xi) d Q_w^-(\xi)\Big)~,
\end{align}
where in the last inequality we applied Lemma \ref{lem-influence-tau}.
Plugging \eqref{eq-bound-f-x} into \eqref{eq-recursion-m_v} now completes the proof of the lemma.
\end{proof}
Observe that in the free boundary case, $Q_v^+(\xi)=Q_v^-(-\xi)$ for any $\xi$. Unfortunately, the presence of the boundary $\tau$ breaks this symmetry,
causing the distributions $Q_v^+$ and $Q_v^-$ to become skewed. Nevertheless, we can still relate these two distributions
through the help of $x_v^*$, as formulated by the following lemma.

\begin{lemma}\label{lem-relation-Q-s}
For $v\in T$, let $Q_v$ be the following distribution over $\{\pm1\}^{\partial\hat{T}}$:
$$Q_v (\xi) = Q_v^\tau(\xi) \deq \mu_v^\tau \left(\sigma_{\partial\hat{T}_v} = \xi_{\partial\hat{T}_v}\right)~.$$
Then for all $\xi \in \{\pm1\}^{\partial\hat{T}}$, we have
$$Q_v^+(\xi) - Q_v^-(\xi) =C^*_v \Big(\tanh\frac{x_v^\xi}{2} - \tanh \frac{x_v^*}{2}\Big) Q_v(\xi)~,$$
where $C^*_v = 2\cosh^2(x_v^*/2)$.
\end{lemma}
\begin{proof}
It is clear from the definitions of $Q_v^+$, $Q_v^-$ and $Q_v$ that
\begin{align*}
Q_v^+(\xi) &= \frac{Q_v (\xi)\mu_v^\tau (\sigma(v) = 1 \given \xi)}{\mu_v^\tau(\sigma(v) = 1)}=\frac{1+\tanh(x_v^\xi/2)}{1+\tanh(x_v^*/2)}Q_v(\xi)~,\\
Q_v^-(\xi) &= \frac{Q_v (\xi)\mu_v^\tau (\sigma(v) = -1 \given \xi)}{\mu_v^\tau(\sigma(v) = -1)}=\frac{1-\tanh(x_v^\xi/2)}{1-\tanh(x_v^*/2)} Q_v(\xi)~.
\end{align*}
Hence, a straightforward calculation gives that
\begin{align*}
Q_v^+(\xi) - Q_v^-(\xi)
&= \frac{2\big(\tanh(x_v^\xi/2) - \tanh(x_v^*/2)\big)}{(1+\tanh(x_v^*/2))(1-\tanh(x_v^*/2))}Q_v(\xi)\\
&=2 \cosh^2\Big(\frac{x_v^*}{2}\Big)\Big(\tanh\frac{x_v^\xi}{2} - \tanh\frac{x_v^*}{2}\Big)Q_v(\xi)~,
\end{align*}
as required.
\end{proof}
The following technical lemma will allow us to combine Lemmas \ref{lem-recursion-m-f}, \ref{lem-relation-Q-s} and obtain an upper bound on $m_v$ in terms of $\{m_w: (v,w)\in \hat{T}\}$. Note that the constant $\kappa$ mentioned next is in fact the absolute constant $\kappa$ in the statement of Theorem \ref{thm-log-likelihood-decay}.

\begin{lemma}\label{lem-f-analysis}
Let $f$ be defined as in \eqref{eq-def-f} for some $0 < \theta \leq \frac34$. Then
\begin{align}\label{eq-f-delta}
  \left|f(x)-f(y)\right| \leq 2f(|x-y|/2)\mbox{ for any $x,y\in\R$}~,
\end{align}
and there exists a universal constant $\kappa > \frac1{100} $ such that for any two constants
$C_1,C_2 \geq 1$ with $C_2 \geq 1 + \left(\frac12C_1-1\right)(1-\theta^2)$ and
any $\delta > 0$ we have
\begin{align}\label{eq-deltaf-times-deltax-bound}
f(\delta)\left(1+ 4 \kappa (1-\theta) C_1 \delta\tanh(\delta/2)\right)
\leq C_2 \theta \delta~.
\end{align}
\end{lemma}
\begin{proof} We first show \eqref{eq-f-delta}. 
Put $\delta = |x-y|$, and define
$$h(\delta) = \sup_{t} \left|f(t+\delta)-f(t)\right|~.$$
We claim that
\begin{equation}\label{eq-h-max}h(\delta) = f(\delta/2) - f(-\delta/2)= 2 f(\delta/2)~.\end{equation}
The second equality follows from the fact that $f(x)$ is an odd function. To establish the first equality, a straightforward calculation gives that
$$f'(x) = \frac{\theta}{1 + (1-\theta^2)\sinh^2(x/2)}~,$$
and it follows that $f'(x)$ is an even non-negative function which is decreasing in $x \geq 0$. The following simple
claim therefore immediately implies \eqref{eq-h-max}:
\begin{claim}
  \label{clm-even-decreasing-derivative}
Let $g(t)\geq 0 $ be an even function that is decreasing in $t\geq 0$. Then $G(t)=\int_0^t g(x)dx$ has
$G(t+\delta) - G(t) \leq 2G(\delta/2)$ for any $t$ and $\delta >0$.
\end{claim}
\begin{proof}
Fix $\delta > 0$ and define $F(t)$ as follows:
$$F (t) = G(t+\delta) - G(t)~.$$
We therefore have $F'(t) = g(t + \delta) - g(t)$. Noticing that
$$\left\{\begin{array}
  {ll}|t + \delta| \geq |t| &\mbox{ if $t \geq -\frac{\delta}{2}$}\\
  |t + \delta| \leq |t| & \mbox{ otherwise}
\end{array}\right.~,$$ the assumption on $g$ now gives that $F'(t) \leq 0$ while $t \geq -\frac{\delta}{2}$ and otherwise $F'(t) \geq 0$. Altogether, we deduce that
$$F(t) \leq F\big(-\mbox{$\frac{\delta}{2}$}\big) =G\big(\mbox{$\frac{\delta}{2}$}\big) - G\big(-\mbox{$\frac{\delta}{2}$}\big)= 2 G\big(\mbox{$\frac{\delta}{2}$}\big)~,$$
as required.
\end{proof}

It remains to show that \eqref{eq-deltaf-times-deltax-bound} holds for some $\kappa=\kappa(\theta_0) > 0$.
Clearly, it suffices to establish this statement for $C_2 = \left[1 + \left(\frac12C_1-1\right)(1-\theta^2)\right] \vee 1$ and any $C_1 \geq 1$.
Rearranging \eqref{eq-deltaf-times-deltax-bound}, we are interested in a lower bound on
\begin{equation}\label{eq-kappa-equation}
\inf_{\theta\leq\theta_0\;,\;t > 0\;,\;C_1 \geq 1} \frac{ \left[\Big(1+\left(\frac{C_1}2 - 1\right)(1-\theta^2)\Big) \vee 1\right] \theta t - f(t) } { 4  C_1 (1-\theta) t f(t)\tanh(\frac{t}2)}~.
\end{equation}
First, consider the case $1 \leq C_1 < 2$. We then have $C_2=1$, and the expression being minimized in \eqref{eq-kappa-equation} takes the form:
$$
\frac{ \theta t - f(t) } { 4 C_1 (1-\theta) t f(t)\tanh(\frac{t}2)}
> \frac{ \theta t - f(t) } { 8 (1-\theta) t f(t)\tanh(\frac{t}2)}\deq\Gamma(t,\theta)~,
$$
where the inequality is by our assumption that $C_1 < 2$. We therefore have that $\inf_{\theta\leq\theta_0\;,\;t > 0} \Gamma(t,\theta)$
minimizes \eqref{eq-kappa-equation} for $C_1 < 2$, and will next show that this is also the case for $C_1 \geq 2$
under a certain condition. Indeed, letting
$$ g(t,\theta,C_1) \deq \frac{ \left[1+\left(\frac{C_1}2 - 1\right)(1-\theta^2)\right] \theta t - f(t) } { 4 C_1 (1-\theta) t f(t)\tanh(t/2)}~,$$
it is easy to verify that the following holds:
\begin{align*}
  \frac{\partial g}{\partial C_1} = \frac{ f(t)-\theta^3 t}{4 C_1^2 (1-\theta)t f(t) \tanh(\frac{t}2)}~,
\end{align*}
hence $g$ is increasing in $C_1$ for every $\theta,t$ such that  $f(t) > \theta^3 t$. Therefore,
$$ g(t,\theta,C_1) \geq g(t,\theta,2) = \Gamma(t,\theta)~\mbox{for all $t,\theta$ such that $f(t)>\theta^3 t$}~.$$
Before analyzing $\Gamma(t,\theta)$, we will treat the values of $\theta,t$ such that $f(t) \leq \theta^3 t$. Assume that
 the case, and notice that the numerator of $g$ then satisfies
\begin{align*}
 \Big[1+\Big(\frac{C_1}2 - 1\Big)&(1-\theta^2)\Big] \theta t - f(t) \\
&\geq  \left[1+\left(\frac{C_1}2 - 1\right)(1-\theta^2)-\theta^2 \right] \theta t = \theta (1 - \theta^2) t \frac{C_1}2~,
\end{align*}
and thereby the dependency on $C_1$ vanishes:
$$ g(t,\theta,C_1) \geq \frac{ \theta (1 - \theta^2)t/2 } { 4 (1-\theta)t f(t)\tanh(t/2)} = \frac{ \theta (1 + \theta) } { 8  f(t)\tanh(t/2)}~.$$
Since both $\tanh(t/2)$ and $f(t)$ are monotone increasing in $t$ and are bounded from above by
$1$ and $\log\left(\frac{1+\theta}{1-\theta}\right)$ respectively, we get
\begin{equation}
  \label{eq-g(t,theta,C1)-inequality}
  g(t,\theta,C_1) \geq \frac{ \theta (1 + \theta)  } { 8   \log\left(\frac{1+\theta}{1-\theta}\right) }
\geq \frac{ \theta (1+\theta) } { 8  \frac{2\theta}{1-\theta} } = \frac{1-\theta^2}{16} > \frac{1}{40}~,
\end{equation}
where the second inequality is by the fact that $\log(1+x) \leq x$ for any $x>0$, and the last inequality follows by the assumption $\theta \leq \frac34$~.

It thus remains to establish a uniform lower bound on $\Gamma(t,\theta)$. In what follows,
our choice of constants was a compromise between simplicity and the quality of the lower bound,
and we note that one can easily choose constants that are slightly more optimal.

Assume first that $\theta \geq \theta_0 \geq 0$ for some $\theta_0$ to be defined later.
Notice that
$$\tilde{f}(t,\theta)\deq \frac1\theta f(t,\theta)= 2\sum_{i=0}^\infty \frac{\tanh^{2i+1}(t/2)}{2i+1}\theta^{2i}~,$$
 and so $\tilde{f}(t,\theta)$ is strictly increasing in $\theta$ for any $t > 0$. Since
\begin{align*} \Gamma(t,\theta) &= \frac{ \theta t - f(t) } { 8 (1-\theta) t f(t)\tanh(\frac{t}2)}\\
&\geq
\frac{ \theta t - f(t) } { 8 (1-\theta_0) t f(t)\tanh(\frac{t}2)} = \frac{ t - \tilde{f}(t) } { 8 (1-\theta_0) t \tilde{f}(t)\tanh(\frac{t}2)}~,\end{align*}
we have that $\Gamma$ is monotone decreasing in $\theta$ for any such $t$, and therefore
 $\Gamma(t,\theta) \geq \frac1{8(1-\theta_0)}\tilde{\Gamma}(t)$, where $\tilde{\Gamma}$ is defined as follows:
 \begin{equation}
   \label{eq-tilde-Gamma-def}
   \tilde{\Gamma}(t) \deq \frac{ \theta t - f(t,\theta) } { t f(t,\theta)\tanh(\frac{t}2)}\quad\mbox{ with respect to $\theta=\frac34$}~.
 \end{equation}
Recall that the Taylor expansion of $f(t,\theta)$ around $0$ is $\theta t - \frac{\theta(1-\theta^2)}{12} t^3 + O(t^5)$.
It is easy to verify that for $\theta=\frac34$ this function satisfies
$$ f(t,\theta) \leq \theta t - \frac{(\theta t)^3}{20}\quad\mbox{ for $\theta=\frac34$ and any $0 < t \leq 3$}~.$$
Adding the fact that $\tanh(x)\leq x$ for all $x\geq 0$, we immediately obtain that
$$\tilde{\Gamma} \geq \frac{\theta^3 t^3}{20 t (\theta t) (t/2)} =\frac{\theta^2}{10} > \frac{1}{20}~\mbox{ for all $0 < t \leq 3$}~.$$
On the other hand, for $t \geq 3$ we can use the uniform upper bounds of $1$ and $\log(\frac{1+\theta}{1-\theta})$ for $\tanh(t/2)$
and $f(t)$ respectively, and gain that
$$\tilde{\Gamma} \geq \frac{\theta t - \log(\frac{1+\theta}{1-\theta})}{t \log(\frac{1+\theta}{1-\theta})} = \frac{\theta}{\log(\frac{1+\theta}{1-\theta})} - \frac1t \geq \frac{1}{20}~\mbox{ for all $t \geq 3$}~.$$
Altogether, as $\Gamma \geq \frac1{8(1-\theta_0)}\tilde{\Gamma}$, we can conclude that $\Gamma \geq \left[160(1-\theta_0)\right]^{-1}$.

Note that the trivial choice of $\theta_0=0$
already provides a uniform lower bound of $\frac1{160}$ for $\Gamma$ (and hence also for $\kappa$, as the lower bound in \eqref{eq-g(t,theta,C1)-inequality} is only larger).
However, this bound can be improved by choosing another $\theta_0$ and treating the case $0 < \theta \leq \theta_0$ separately.
To demonstrate this, take for instance $\theta_0 = \frac12$. Since the above analysis gave that $\tilde{\Gamma} \geq \frac{1}{20}$
whenever $\theta \leq \frac34$, it follows that
$$\Gamma \geq \frac1{160(1-\theta_0)} = \frac1{80}\quad\mbox{ for all $\frac12 \leq \theta \leq \frac34$}~.$$
For $\theta \leq \theta_0$, we essentially repeat this analysis of $\tilde{\Gamma}$, only this time the respective value
of $\theta$ (that is, the maximum value it can attain) is $\frac12$. One can thus verify that in that case,
$$ f(t,\theta) \leq \theta t - \frac{(\theta t)^3}{6}\quad\mbox{ for $\theta=\frac12$ and any $0 < t \leq 2.7$}~,$$
and the above argument then shows that
$$\tilde{\Gamma} \geq \frac{\theta^2 }{3} =\frac{1}{12}~\mbox{ for all $0 < t \leq 2.7$}~.$$
On the other hand,
$$\tilde{\Gamma} \geq  \frac{\theta}{\log(\frac{1+\theta}{1-\theta})} - \frac1t \geq \frac{1}{12}~\mbox{ for all $t \geq 2.7$}~,$$
thus for $\theta=\frac12$ we have $\tilde{\Gamma} \geq \frac1{12}$ for all $t > 0$. This converts into the lower bound $\Gamma \geq \frac1{96}$, thus completing the proof with a final value of $\kappa=\frac1{96}$.
\end{proof}
\begin{remark}\label{rem-kappa-theta-3/4} Note that the only places where we used the fact that $\theta \leq \frac34$ are the lower bound on $g(t,\theta,C_1)$ in \eqref{eq-g(t,theta,C1)-inequality} and the analysis of $\tilde{\Gamma}$, as defined in \eqref{eq-tilde-Gamma-def}. In both cases, we actually only need to have $\theta \leq \theta_1 $ for some constant $\theta_1 < 1$, whose precise value might affect
the final value of $\kappa$.
\end{remark}
Using the above lemma, we are now ready to obtain the final ingredient required for the proof of the recursion relation of $m_v$, as incorporated in Lemma~\ref{lem-recursion-f-m}. This lemma provides a recursive bound on a quantity that resembles $m_v$, where
instead of integrating over $x_v^\xi$, we integrate over $f(x_v^\xi)$.
\begin{lemma}\label{lem-recursion-f-m}
Let $f$ and $D_v^*$ be as in \eqref{eq-def-f} and Lemma \ref{lem-influence-tau} respectively. There exists a universal constant $\kappa > \frac1{100}$  so that for $K = \frac14 (1-\theta)\kappa$ we have
$$\int f(x_v^\xi) d Q_v^+(\xi) - \int f(x_v^\xi) d Q_v^-(\xi) \leq  \frac{\theta m_v}{D_v^*(1 + K m_v)}~.$$
\end{lemma}
\begin{proof}
Clearly,
\begin{align*}
\int f(x_v^\xi) d Q_v^+(\xi) &- \int f(x_v^\xi) d Q_v^-(\xi)\\
&= \int (f(x_v^\xi) - f(x_v^*)) d Q_v^+(\xi) - \int (f(x_v^\xi) - f(x_v^*)) d Q_v^-(\xi)~.
\end{align*}
Applying Lemma \ref{lem-relation-Q-s}, we then obtain that
\begin{align*}
\int f(x_v^\xi) d Q_v^+(\xi) &- \int f(x_v^\xi) d Q_v^-(\xi)
\\ &=  C_v^* \int (f(x_v^\xi) - f(x_v^*))\Big(\tanh\frac{x_v^\xi}{2} - \tanh\frac{x_v^*}{2}\Big) d Q_v(\xi)~,
\end{align*}
and similarly,
$$m_v  = C_v^* \int (x_v^\xi - x_v^*)\Big(\tanh\frac{x_v^\xi}{2} - \tanh \frac{x_v^*}{2}\Big) d Q_v(\xi)~.$$
Let
\begin{align*}
F(x)& = \left(f(x)-f(x_v^*)\right)\left(\tanh(x/2) - \tanh(x_v^*/2)\right)~,\\
G(x) &= \left(x - x_v^*\right)\left(\tanh(x/2) - \tanh(x_v^*/2)\right)~,
\end{align*}
and define $\Lambda$ to be the probability measure on $\R$ as:
$$ \Lambda(x) \deq Q_v\left(\{ \xi: x_v^\xi = x\}\right)~.$$
According to this definition, we have
\begin{align*}\int F(x_v^\xi) d Q_v(\xi) =\int F(x) d \Lambda~,~\mbox{ and }\int G(x_v^\xi) d Q_v(\xi) =\int G(x) d \Lambda~,
\end{align*}
and thus, by the above arguments,
\begin{align}
&\int f(x_v^\xi) d Q_v^+(\xi) - \int f(x_v^\xi) d Q_v^-(\xi) = C_v^* \int F(x) d\Lambda~,\nonumber\\
 & m_v = C_v^* \int G(x) d\Lambda~.\label{eq-mv-lambda}
\end{align}
Furthermore, notice that by \eqref{eq-f-delta} and the fact that $f$ is odd and increasing for $x\geq0$,
\begin{equation*}
 F(x) \leq 2f\Big(\frac{x-x_v^*}2\Big)\Big(\tanh\frac{x}2 - \tanh\frac{x_v^*}2\Big)~.
\end{equation*}
and so
\begin{align}
  \label{eq-first-F-bound}
\int& f(x_v^\xi) d Q_v^+(\xi) - \int f(x_v^\xi) d Q_v^-(\xi) \nonumber \\
&\leq 2C_v^* \int f\Big(\frac{x-x_v^*}2\Big)\Big(\tanh\frac{x}2 - \tanh\frac{x_v^*}2\Big)d\Lambda~.
\end{align}

In our next argument, we will estimate $\int G(x)d\Lambda$ and $\int G(x)d\Lambda$ according to the behavior of $F$ and $G$ about $x_v^*$. Assume that $x_v^* \geq 0$, and note that, although the case of $x_v^* \leq 0$ can be treated similarly, we claim that this assumption does not lose generality. Indeed,  if $x_v^* < 0$, one can consider the boundary condition of $\tau' = -\tau$, which would give the following by symmetry:
$${x_v^*}' = -x_v^*~,~X'_v(-\xi) = -x_v^\xi(\xi)~,~ Q'_v(-\xi)=Q_v(\xi)~.$$
Therefore, as $f(\cdot)$ and $\tanh(\cdot)$ are both odd functions, we have that $\int F(x)d\Lambda$ and $\int G(x)d\Lambda$ will not change under the modified boundary condition, and yet ${x_v^*}' \geq 0$ as required.

Define
$$I^- \deq (-\infty,x_v^*]\quad,\quad I^+ \deq [x_v^*,\infty)~.$$
First, consider the case where for either $I = I^+$ or $I=I^-$ we have
\begin{equation}
  \label{eq-FG-case-1}
  \left\{\begin{array}{l}
\int_{I} F(x) d \Lambda \geq \frac12\int F(x) d \Lambda~,\\
\int_{I} G(x) d \Lambda \geq \frac12\int F(x) d \Lambda~.
\end{array}\right.
\end{equation}
In this case, the following holds
\begin{align*}
\left(\int F(x) d \Lambda\right)\left(\int G(x) d \Lambda\right)
&\leq 4 \left(\int_{I} F(x) d\Lambda\right) \left(\int_{I} G(x) d\Lambda\right)\\
&\leq 4 \int_{I} F(x) G(x) d \Lambda \leq 4 \int F(x) G(x) d \Lambda~,
\end{align*}
where in the second line we applied the FKG-inequality, using the fact that both $F$ and $G$ are decreasing in $I^-$ and increasing in $I^+$. The last inequality followed from the fact that $F$ and $G$ are always non-negative. Note that
$$ \int F(x)G(x)d\Lambda = \int \left(f(x)-f(x_v^*)\right)(x-x_v^*)\Big(\tanh\frac{x}2-\tanh\frac{x_v^*}2\Big)^2 d\Lambda~,$$
and recall that Claim \ref{clm-even-decreasing-derivative} applied onto $\tanh(x)$ (which indeed has an even non-negative derivative $\cosh^{-2}(x)$ that is decreasing in $x \geq 0$) gives
$$\tanh \frac{x}2-\tanh\frac{y}2 \leq
2\tanh \Big(\frac{x-y}4\Big)~\mbox{ for any $x > y$}~. $$
Noticing that each of the factors comprising $F(x)G(x)$ has the same sign as that of $(x-x_v^*)$, and combining this
with \eqref{eq-f-delta}, it thus follows that
\begin{align}  \label{eq-intF-intG-case-1} &\left(\int F(x) d \Lambda\right)\left(\int G(x) d \Lambda\right) \nonumber\\
&\leq 16\int f\Big(\frac{x-x_v^*}2\Big)(x-x_v^*)\Big(\tanh\frac{x}2-\tanh\frac{x_v^*}2\Big)
\tanh\Big(\frac{x-x_v^*}{4}\Big) d\Lambda~. \end{align}


Second, consider the case where for $I^+$ and $I^-$ as above, we have
\begin{equation}
  \label{eq-FG-case-2}
\left\{\begin{array}{l}
\int_{I^+} F(x) d \Lambda \geq \frac12\int F(x) d \Lambda~,\\
\int_{I^-} G(x) d \Lambda \geq \frac12\int G(x) d \Lambda~.
\end{array}\right.
\end{equation}
The following definitions of $\tilde{F}$ and $\tilde{G}$ thus capture a significant contribution of $F$ and $G$ to $\int F d\Lambda$
and $\int Gd\Lambda$ respectively:
\begin{align}\label{eq-tildeF-tildeG-def}
\left\{\begin{array}{l}
\tilde{F}(s) \deq F(x_v^* + s)\\
\tilde{G}(s) \deq G(x_v^* - s)
\end{array}\right.~\mbox{ for any $s \geq 0$}~,
\end{align}
By further defining the probability measure $\tilde{\Lambda}$ on $[0,\infty)$ to be
\begin{equation}\label{eq-def-tilde-Lambda}\tilde{\Lambda}(s) \deq  \Lambda(x_v^*-s)\one_{\{s\neq0\}}+\Lambda(x_v^*+s)~\mbox{ for any $s \geq 0$}~,\end{equation}
  we obtain that
\begin{align*}
\int F(x) d\Lambda&\leq 2\int_{I^+} F(x) d\Lambda \leq 2\int_{0}^\infty\tilde{F}(x) d\tilde{\Lambda}~,\\
\int G(x) d\Lambda&\leq 2\int_{I^-} G(x) d\Lambda \leq 2\int_{0}^\infty\tilde{G}(x) d\tilde{\Lambda}~.
\end{align*}
With both $\tilde{F}$ and $\tilde{G}$ being monotone increasing on $[0,\infty)$, applying
the FKG-inequality with respect to $\tilde{\Lambda}$ now gives
\begin{align*}
&\left(\int F(x) d \Lambda\right)\left(\int G(x) d \Lambda\right) \leq 4 \int_0^\infty \tilde{F}(x) \tilde{G}(x) d \tilde{\Lambda}\nonumber\\
&= 4\int_0^{\infty} \left(f(x_v^*+s)-f(x_v^*)\right)\left(\tanh\frac{x_v^*+s}2-\tanh\frac{x_v^*}2\right)\nonumber\\
&\qquad\quad\cdot (-s)\left(\tanh\frac{x_v^*-s}2-\tanh\frac{x_v^*}2\right) d\tilde{\Lambda}~.
\end{align*}
Returning to the measure $\Lambda$, the last expression takes the form
\begin{align*}
 4&\int_{I^+} \left(f(x)-f(x_v^*)\right)\left(\tanh\frac{x}2-\tanh\frac{x_v^*}2\right)\\
&\quad\cdot(x-x_v^*)\left(\tanh\frac{x_v^*}2-\tanh\frac{2x_v^*-x}2\right) d\Lambda\\
+4&\int_{I^-} \left(f(2x_v^*-x)-f(x_v^*)\right)\left(\tanh\frac{2x_v^*-x}2-\tanh\frac{x_v^*}2\right)\\
&\quad\cdot(x-x_v^*)\left(\tanh\frac{x}2-\tanh\frac{x_v^*}2\right) d\Lambda~.\\
\end{align*}
We now apply \eqref{eq-f-delta} and Claim \ref{clm-even-decreasing-derivative} (while leaving the term $(\tanh\frac{x}2-\tanh\frac{x_v^*}2)$ unchanged in both integrals) to obtain that
\begin{align*}
&\left(\int F(x) d \Lambda\right)\left(\int G(x) d \Lambda\right) \nonumber\\
&\leq 16\int f\Big(\frac{x-x_v^*}2\Big)\left(\tanh\frac{x}2-\tanh\frac{x_v^*}2\right)(x-x_v^*)
\tanh\Big(\frac{x-x_v^*}4\Big) d\Lambda~.
\end{align*}
That is, we have obtained the same bound as in \eqref{eq-intF-intG-case-1}.

It remains to deal with the third case where for $I^+$ and $I^-$ as above,
\begin{equation}\label{eq-FG-case-3}
\left\{\begin{array}{l}
\int_{I^-} F(x) d \Lambda \geq \frac12\int F(x) d \Lambda~,\\
\int_{I^+} G(x) d \Lambda \geq \frac12\int G(x) d \Lambda~.
\end{array}\right.\end{equation}
In this case, we modify the definition \eqref{eq-tildeF-tildeG-def} of $\tilde{F}$ and $\tilde{G}$ appropriately:
\begin{align*}\left\{\begin{array}{l}
\tilde{F}(s) \deq F(x_v^* + s)\\
\tilde{G}(s) \deq G(x_v^* - s)
\end{array}\right.~\mbox{ for any $s \geq 0$}~,
\end{align*}
and let $\tilde{\Lambda}$ remain the same, as given in \eqref{eq-def-tilde-Lambda}. It then follows that
\begin{align*}
\int F(x) d\Lambda&\leq 2\int_{I^-} F(x) d\Lambda \leq 2\int_{0}^\infty\tilde{F}(x) d\tilde{\Lambda}~,\\
\int G(x) d\Lambda&\leq 2\int_{I^+} G(x) d\Lambda \leq 2\int_{0}^\infty\tilde{G}(x) d\tilde{\Lambda}~,
\end{align*}
with $\tilde{F}$ and $\tilde{G}$ monotone increasing on $[0,\infty)$. By the FKG-inequality,
\begin{align}
&\left(\int F(x) d \Lambda\right)\left(\int G(x) d \Lambda\right) \leq 4 \int_0^\infty \tilde{F}(x) \tilde{G}(x) d \tilde{\Lambda}\nonumber\\
&= 4\int_0^{\infty} \left(f(x_v^*-s)-f(x_v^*)\right)\left(\tanh\frac{x_v^*-s}2-\tanh\frac{x_v^*}2\right)\nonumber\\
&\qquad\qquad\cdot s\left(\tanh\frac{x_v^*+s}2-\tanh\frac{x_v^*}2\right) d\tilde{\Lambda}~.\label{eq-intF-intG-prelim-bound-case-3}
\end{align}
As before, we now switch back to $\Lambda$ and infer from \eqref{eq-f-delta} and Claim \ref{clm-even-decreasing-derivative} that
\begin{align*}
&\left(\int F(x) d \Lambda\right)\left(\int G(x) d \Lambda\right)\\
&\leq 16\int f\Big(\frac{x-x_v^*}2\Big)\left(\tanh\frac{x}2-\tanh\frac{x_v^*}2\right)(x-x_v^*)
\tanh\Big(\frac{x-x_v^*}4\Big) d\Lambda~,
\end{align*}
that is, \eqref{eq-intF-intG-case-1} holds for each of the 3 possible cases \eqref{eq-FG-case-1}, \eqref{eq-FG-case-2} and \eqref{eq-FG-case-3}.

Altogether, this implies that
\begin{align*}
&\left(\int f(x_v^\xi) d Q_v^+(\xi) - \int f(x_v^\xi) d Q_v^-(\xi) \right) m_v \\
&= (C_v^*)^2\left(\int F(x)d\Lambda\right)\left(\int G(x)d\Lambda\right)\\
& \leq 16{C_v^*}^2\int f\Big(\frac{x-x_v^*}2\Big)\left(\tanh\frac{x}2-\tanh\frac{x_v^*}2\right)(x-x_v^*)
\tanh\Big(\frac{x-x_v^*}4\Big) d\Lambda~.
\end{align*}
Therefore, recalling \eqref{eq-first-F-bound} and choosing $K = \frac14(1-\theta)\kappa$, where $\kappa$ is as given
in Lemma \ref{lem-f-analysis}, we have
\begin{align*}
&\Big(\int f(x_v^\xi) d Q_v^+(\xi) - \int f(x_v^\xi) d Q_v^-(\xi)\Big)(1 + K m_v)\\
&\leq 2C_v^* \int f\Big(\frac{x-x_v^*}2\Big)\Big[1+ 4\kappa(1-\theta)  C_v^* \frac{x-x_v^*}2
\tanh\frac{x-x_v^*}4\Big] \\
&\qquad\qquad\cdot\Big(\tanh\frac{x}2 - \tanh\frac{x_v^*}2\Big)d\Lambda \\
%
&\leq 2C_v^* \int (1/D_v^*) \theta \frac{x-x_v^*}2\Big(\tanh\frac{x}2 - \tanh\frac{x_v^*}2\Big)d\Lambda
= \theta \frac{C_v^*}{D_v^*} \int G(x)d\Lambda~.
\end{align*}
where the inequality in the last line is by Lemma \ref{lem-f-analysis} for $\delta = |x-x_v^*|/2$ (the case
 $x < x_v^*$ follows once again from the fact that $f$ is odd) and a choice of
 $C_1 = C_v^* = 2\cosh^2(x_v^*/2) \geq 2$
and $C_2 = (1/D_v^*)$ (recall that, by definition, $1/D_v^* = 1 + (\frac12 C_v^* -1 )(1-\theta^2)\geq 1$, satisfying the requirements of the lemma).
Therefore, \eqref{eq-mv-lambda} now implies that
\begin{align*}
&\int f(x_v^\xi) d Q_v^+(\xi) - \int f(x_v^\xi) d Q_v^-(\xi)\leq \frac{\theta m_v}{D_v^*(1 + K m_v) }~,
\end{align*}
as required.
\end{proof}
Combining Lemmas \ref{lem-recursion-m-f} and \ref{lem-recursion-f-m}, we deduce that there exists a universal constant $\kappa > 0$ such that
\begin{equation}\label{eq-mv-recursion}m_v \leq \sum_{w:(v,w)\in\hat{T}} \frac{\theta^2 m_w}{ 1 + \frac14 \kappa(1-\theta) m_w}~.\end{equation}
The proof will now follow from a theorem of \cite{PP}, that links a function on the vertices of a tree $T$ with its $L^2$-capacity according to certain resistances.
\begin{theorem}[\cite{PP}*{Theorem 3.2} (reformulated)]\label{thm-recursion-cap}
Let $T$ be a finite tree, and suppose that there exists some $K > 0$ and positive constants $\{a_v: v\in T\}$ such that for every $v\in T$ and $x \geq 0$,
$$ g_v(x) \leq a_v x \big/(1+ K x) ~.$$
Then any solution to the system $x_v = \sum_{w:(v,w)\in T} g_w (x_w) $
satisfies
$$x_\rho \leq \mathrm{cap}_2 (T) \big/ K~,$$
where the resistances are given by $R_{(u,v)} = \prod_{(x,y)\in\mathcal{P}(\rho,v)} a_y^{-1}$, with $\mathcal{P}(\rho,v)$ denoting
the simple path between $\rho$ and $v$.
\end{theorem}
Together with inequality \eqref{eq-mv-recursion}, the above theorem immediately gives
$$m_\rho \leq \frac{\mathrm{cap}_2 (\hat{T})}{\kappa (1-\theta)/4}~,$$
completing the proof of Theorem \ref{thm-log-likelihood-decay}. \qed

\subsection*{Proof of Proposition \ref{prop-delta}}
In order to obtain the required result from Theorem~\ref{thm-log-likelihood-decay},
recall the definition of $x_v^\xi$ for $v\in T$, according to which we can write
$$\hat{\mu}^\xi(\sigma(\rho) = 1) = \left(1+ \tanh(x_\rho^\xi/2 + h)\right)/2~,$$
where $h$ is the mentioned external field at the root $\rho$.
By monotone coupling, we can construct a probability measure $Q_c$ on the space $\{(\xi, \xi'): \xi \geq \xi'\}$ such that the two marginal distributions correspond to $Q_\rho^+$ and $Q_\rho^-$ respectively. It therefore follows that
\begin{align*}
\Delta &= \int \Big(\hat{\mu}^\xi(\sigma(\rho) = 1) - \hat{\mu}^{\xi'}(\sigma(\rho) = 1)\Big) d Q_c\\
&=\frac12\int  \Big(\tanh(x_\rho^\xi/2 +h) - \tanh(x_\rho^{\xi'} /2 +h) \Big)d Q_c\\
&\leq \frac12\int \frac{x_\rho^\xi - x_\rho^{\xi'}}{2} dQ_c = \frac{1}{4} m_\rho
\leq \frac{\mathrm{cap}_2 (\hat{T})}{\kappa (1-\theta)}~,
\end{align*}
where the last inequality follows from Theorem \ref{thm-log-likelihood-decay} using the same value of $\kappa \geq \frac1{100}$.
This completes the proof.
\qed

\section{Upper bound on the inverse-gap and mixing time}\label{sec:thm-1}
This section is devoted to the proof of the main theorem, Theorem \ref{thm-1}, from which it
follows that the mixing time of the continuous-time Glauber dynamics for the Ising model on a $b$-ary tree (with any boundary condition)
is poly-logarithmic in the tree size.

Recalling the log-Sobolev results described in Section \ref{sec:prelim}, it suffices to show an upper bound of $O(n \log^M n)$ on inverse-gap of the discrete-time chain (equivalently, a lower bound on its $\gap$), which would then imply an upper bound of $O(n \log^{M+2} n)$ for the $L^2$ mixing-time (and hence also for the total-variation mixing-time).

The proof comprises several elements, and notably, uses a block dynamics in order to obtain the required upper bound inductively.
Namely, we partition a tree on $n$ vertices to blocks of size roughly $n^{1-\alpha}$ each, for some small $\alpha>0$, and use
an induction hypothesis that treats the worst case boundary condition. The main effort is then to establish a lower bound
on the spectral-gap of the block dynamics (as opposed to each of its individual blocks). This is achieved by Theorem \ref{thm-gap-block-dynamics} (stated later), whose proof hinges on the spatial-mixing result of Section \ref{sec:spatial-mixing}, combined with the Markov chain decomposition method.

Throughout this section, let $b \geq 2$ be fixed, denote by $\beta_c=\arctanh({1/\sqrt{b}})$ the critical inverse-temperature and let $\theta = \tanh \beta_c$.

\subsection{Block dynamics for the tree}
In what follows, we describe our choice of blocks for the above mentioned block dynamics. Let $h$ denote the height of our $b$-ary tree (that is, there are $b^{h}$ leaves in the tree),
and define
\begin{equation}
  \label{eq-lup-ldown-def}
\lup \deq \alpha h \quad, \quad \ldown \deq h - \lup~,
\end{equation}
where $0 < \alpha < \frac12$ is some (small) constant to be selected later.

For any $v \in T$, let $B(v,k)$ be the subtree of height $k-1$ rooted at $v$, that is, $B(v, k)$ consists of $k$ levels (except when $v$ is less than $k$ levels away from the bottom of $T$). We further let $H_k$ denote the $k$-th level of the tree $T$, that according
to this notation contains $b^k$ vertices.

Next, define the set of blocks $\mathcal{B}$ as:
\begin{equation}
  \label{eq-block-dynamics-def}
\mathcal{B} \deq \left\{B(v, \ldown)~:~ v \in H_\lup \cup \{\rho\}\right\}~\quad\mbox{for $\lup,\ldown$ as above.}
\end{equation}
That is, each block is a $b$-ary tree with $r$ levels, where one of these blocks is rooted at $\rho$, and will be referred
to as the \emph{distinguished block}, whereas the others are rooted at the vertices of $H_\lup$.

\begin{figure}
\centering \includegraphics[width=3in]{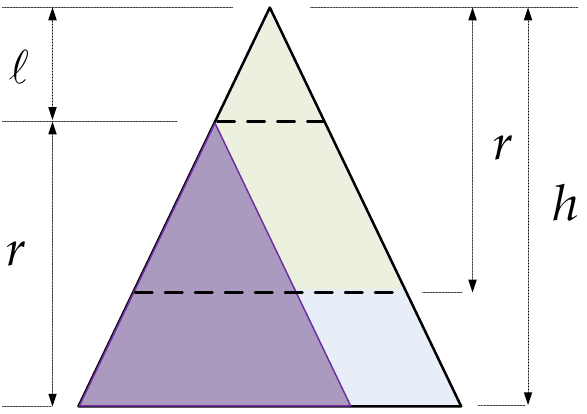}
\caption{Block dynamics for the Ising model on the tree: illustration shows
the distinguished block $B(\rho,r)$ as well as a representative block of the form $B(v,r)$ for $v\in H_\lup$.}
\label{fig:block-dynamics}
\end{figure}

The following theorem establishes a lower bound on the spectral gap of the above-specified block dynamics (with blocks $\mathcal{B}$).
\begin{theorem}\label{thm-gap-block-dynamics}
Consider the Ising model on the $b$-ary tree at the critical inverse-temperature $\beta_c$ and with an arbitrary boundary $\tau$.
Let $\gap_\mathcal{B}^\tau$ be the spectral gap of the corresponding block dynamics with blocks $\mathcal{B}$ as in \eqref{eq-block-dynamics-def}.
The following then holds:
$$\gap_\mathcal{B}^\tau \geq \frac{1}{4(b^{\lup} + 1)} \bigg(1 - \frac{\alpha}{\kappa(1-\theta)(1-2\alpha)}\bigg)~,$$
where $\kappa>0$ is the absolute constant given in Theorem \ref{thm-log-likelihood-decay}.
\end{theorem}
Given the above theorem, we can now derive a proof for the main result.
\subsection{Proof of Theorem \ref{thm-1}}
By definition, as $b \geq 2$, we have that $$\theta = \tanh \beta_c =  \frac{1}{\sqrt{b}} \leq \frac{1}{\sqrt{2}}~,$$ hence we can readily choose an absolute constant $0<\alpha<1$ such that
$$c(\alpha) \deq \frac{1}{8} \bigg(1 -\frac{\alpha}{\kappa(1-\theta)(1-2\alpha)}\bigg) > 0~.$$
Let $n_h=\sum_{j=0}^{h-1}b^j $ be the number of vertices in a $b$-ary tree of height $h$ excluding its leaves, and let $\gap_h^\tau$ be the spectral gap of the (single-site) discrete-time Glauber dynamics for the Ising model on a $b$-ary tree of height $h$ and boundary $\tau$ (in the special case of a free boundary condition, $n_h$ should instead include the leaves). Further define
$$g_h = n_h \min_{\tau} \gap_h^\tau~.$$
Recalling the definition of $\mathcal{B}$ according to the above choice of $\alpha$, we have that each of its blocks is a tree of height $r =(1-\alpha)h$, and that
$$\sup_{v\in T} \#\{B \in \mathcal{B}~:~ x\in B\} = 2~,$$
as each of the vertices in levels $\lup,\lup+1,\ldots,\ldown$ is covered precisely twice in $\mathcal{B}$, while every other vertex is covered precisely once.

Hence, by Proposition \ref{prop-block-single} and Theorem \ref{thm-gap-block-dynamics}, it now follows that for any
$h \geq 1/\alpha$ (such that our choices of $\lup,\ldown$ in \eqref{eq-lup-ldown-def} are both non-zero) we have
$$g_h \geq \bigg(\frac{1}{4(b^{\lup} + 1)} \bigg(1 - \frac{\alpha}{\kappa(1-\theta)(1-2\alpha)}\bigg) \bigg) g_{\ldown} \cdot\frac12  = c(\alpha) g_{(1-\alpha) h}~.$$
Having established the induction step, we now observe that, as $\alpha$ is constant,
clearly $g_k \geq c'$ holds for any $k \leq 1/\alpha$ and some fixed $c' = c'(\alpha) > 0$. Hence,
$$ g_h \geq c' \big(c(\alpha)\big)^{\log_{1-\alpha}(1/h)} = c' h^{-\log\left(\frac1{c(\alpha)}\right) / \log \left(\frac{1}{1-\alpha}\right)}~,$$
that is, there exists an absolute constant $M$ (affected by our choice of the absolute constants $\kappa,\alpha$) so
that the inverse-gap of the continuous-time dynamics with an arbitrary boundary condition $\tau$ is at most $g_h^{-1} = O(h^M)$, as required.
\qed

\subsection{Proof of Theorem \ref{thm-gap-block-dynamics}}
In order to obtain the desired lower bound on the spectral gap of the block dynamics, we will apply the method of decomposition of Markov chains, described in Subsection \ref{subsec:prelim-decomposition-mc}. To this end, we will partition our configuration according to the spins of the subset
 $$S \deq B(\rho, \lup-1)~.$$
Note that $S$ is strictly contained in the distinguished block $B(\rho,r)$, and does not intersect any other $B \in \mathcal{B}$.
For $\eta \in \{\pm1\}^S$, denote the set of configurations which agree with $\eta$ by
$$\Omega_\eta \deq \{\sigma\in \Omega: \sigma_S = \eta\}~.$$
 Following the definitions in Subsection \ref{subsec:prelim-decomposition-mc}, we can now naturally decompose the block dynamics into a projection chain $\bar{P}$ on $\{\pm1\}^S$ and restriction chains $P_{\eta}$ on $\Omega_\eta$ for each $\eta\in\{\pm1\}^S$. With Theorem \ref{thm-gap-decomposition} in mind, we now need to provide suitable lower bounds on $\bar{\gap}^\tau$ and $\gap_\eta^\tau$, the respective spectral gaps of $\bar{P}$ and $P_{\eta}$  given the boundary condition $\tau$.

We begin with the lower bound on the restriction chain $\gap^\tau_\eta$, formulated in the next lemma.

\begin{lemma}\label{lem-gap-restr}
For any boundary $\tau$ and $\eta\in \{\pm1\}^S$, the spectral gap of the restriction chain $P_\eta$ on the space $\Omega_\eta$ satisfies
$\gap^\tau_\eta \geq 1/(b^{\lup} + 1)$.
\end{lemma}
\begin{proof}
Recall that the restriction chain $P_\eta$ moves from $\sigma\in\Omega_\eta$ to $\sigma'\in\Omega_\eta$ (that is, $\sigma$ and $\sigma'$ both  agree with $\eta$ on $S$) according to the original law of the chain, and remains at $\sigma$ instead of moving to any $\sigma'\notin \Omega_\eta$.
By definition of our block dynamics, this means that with probability $b^\lup/(b^\lup + 1)$ we apply a transition kernel $Q_1$, that selects one of the blocks rooted at $H_\lup$ to be updated according to its usual law (since $S$ and all of these blocks are pairwise disjoint). On the other hand, with probability $1/(b^\lup + 1)$, we apply a transition kernel $Q_2$ that updates the distinguished block, yet only allows updates that keep $S$ unchanged (otherwise, the chain remains in place).

We next claim that the update of the distinguished block can only increase the value of $\gap^\tau_\eta$. To see this, consider the chain $P'_\eta$, in which the distinguished block is never updated; that is, $Q_2$ described above is replaced by the identity. Clearly, since each of the vertices of $T \setminus S$ appears in (precisely) one of the non-distinguished blocks, the stationary distribution of $P'_\eta$ is again $\mu_{\tau;\eta}$, the Gibbs distribution with boundary conditions $\eta$ and $\tau$. Therefore, recalling the Dirichlet form \eqref{eq-dirichlet-form-e}, for any $f$ we clearly have
\begin{align*}
\mathcal{E}_{P'_\eta}(f) &= \frac{1}{2}\sum_{x,y\in\Omega_\eta}\left[f(x)-f(y)\right]^2\mu^{\tau;\eta}(x)P'_\eta(x,y)\\
&\leq \frac{1}{2}\sum_{x,y\in\Omega_\eta}\left[f(x)-f(y)\right]^2\mu^{\tau;\eta}(x)P_\eta(x,y) = \mathcal{E}_{P_\eta}(f)~,
\end{align*}
and thus, by the spectral gap bound in terms of the Dirichlet form \eqref{eq-dirichlet-form},
\begin{equation}
  \label{eq-gap-Peta-P'eta}
\gap( P_\eta) \geq \gap (P'_\eta)~.
\end{equation}
It remains to analyze the chain $P'_\eta$, which is in fact a \emph{product chain}, and as such its eigenvalues can be directly expressed in terms of the eigenvalues of its component chains. This well known fact is stated in the following straightforward claim (cf., e.g., \cite{AF}*{Chapter 4} and \cite{LPW}*{Lemma 12.11}); we include its proof for completeness.
\begin{claim}\label{claim-product-chain}
For $j \in[d]$, let $P_j$ be a transition kernel on $\Omega_j$ with eigenvalues $\Lambda_j$. Let $\nu$ be a probability distribution on $[d]$, and define $P'$, the transition matrix of the product chain of the $P_j$-s on $\Omega' = \Omega_1\times \Omega_2 \times \cdots \times \Omega_d$, by
$$P'\big((x_1,\ldots,x_d),(y_1,\ldots,y_d)\big) = \sum_{j=1}^d \nu(j) P_j(x_j, y_j) \prod_{i: i\neq j} \one_{\{x_i = y_i\}}~.$$
Then $P'$ has eigenvalues $\left\{ \sum_{j=1}^d \nu(j) \lambda_j : \lambda_j\in \Lambda_j\right\}$ (with multiplicities).
\end{claim}
\begin{proof}
Clearly, by induction it suffices to prove the lemma for $d=2$. In this case, it is easy to verify that
the transition kernel $\tilde{P}$ can be written as
$$ \tilde{P} = \nu(1) (P_1 \otimes I_{\Omega_2}) + \nu(2) (I_{\Omega_1} \otimes P_2)~,$$
where $\otimes$ denotes the matrix tensor-product. Thus, by tensor arithmetic,
for any $u,v$, eigenvectors of $P_1,P_2$ with corresponding eigenvalues
$\lambda_1,\lambda_2$ respectively, $(u\otimes v)$ is an eigenvector of $\tilde{P}$ with
a corresponding eigenvalue of $\nu(1) \lambda_1  + \nu(2)\lambda_2$, as required.
\end{proof}
In our setting, first notice that $Q_1$ itself is a product chain, whose components are the $b^\lup$ chains, uniformly selected, updating each of the non-distinguished blocks. By definition, a single block-update replaces the contents of the block with a sample according to the stationary distribution conditioned on its boundary. Therefore, each of the above mentioned component chains has a single eigenvalue of $1$ whereas all its other eigenvalues are $0$.

It thus follows that $P'_\eta$ (a lazy version of $Q_1$) is another product chain, giving $Q_1$ probability $b^\lup/(b^\lup+1)$ and the identity chain probability $1/(b^\lup+1)$. By Claim \ref{claim-product-chain}, we conclude that the possible eigenvalues of $P'_\eta$ are precisely
$$\left\{ \frac{1}{b^\lup+1} +\frac{1}{b^\lup+1}\mbox{$\sum_{j=1}^{b^\lup} \lambda_j$} : \lambda_j\in\{0,1\}\right\}~.$$
In particular, $\gap(P'_\eta) = 1/(b^\lup + 1)$, and \eqref{eq-gap-Peta-P'eta} now completes the proof.
\end{proof}

It remains to provide a bound on $\bar{\gap}^\tau$, the spectral gap of the projection chain in the decomposition of the block dynamics according to $S$. This is the main part of our proof of the lower bound for the spectral gap of the block dynamics, on which the entire proof of Theorem \ref{thm-1} hinges. To obtain this bound, we relate the projection chain to the spatial-mixing properties of the critical Ising model on the tree under various boundary conditions, studied in Section \ref{sec:spatial-mixing}.

\begin{lemma}\label{lem-gap-proj}
For any boundary $\tau$, the spectral gap of the projection chain $\bar{P}$ on the space $\{\pm1\}^S$ satisfies
$$ \bar{\gap}^\tau \geq \frac1{b^\lup +1}\left( 1 - \frac{\alpha}{\kappa(1-\theta)(1-2\alpha)}\right)~,$$
where $\kappa>0$ is the absolute constant given in Proposition \ref{prop-delta}.
\end{lemma}
We prove this lemma by establishing a certain contraction property of the projection chain $\bar{P}$.
Recall that $\bar{P}(\eta,\eta')$, for $\eta,\eta'\in\{\pm1\}^S$, is the probability that completing $\eta$
into a state $\sigma$ according to the stationary distribution (with boundary $\eta$ and $\tau$) and then
applying the block dynamics transition, gives some $\sigma'$ that agrees with $\eta'$ on $S$.

Let $S^* = H_{\lup -1}$ denote the bottom level of $S$, and notice that in the above definition of the transition kernel of $\bar{P}$, the value of the spins in $S \setminus S^*$ do not affect the transition probabilities. Therefore, the projection of
the chain $\bar{P}$ onto $S^*$ is itself a Markov chain, which we denote by $\bar{P}^*$. In fact,
we claim that the eigenvalues of $\bar{P}$ and those of $\bar{P}^*$ are precisely the same (with the exception of
additional $0$-eigenvalues in $\bar{P}$). To see this, first notice that the eigenfunctions of $\bar{P}^*$
can be naturally extended into eigenfunctions of $\bar{P}$ with the same eigenvalues (as $\bar{P}^*$ is a projection of $\bar{P}$).
Furthermore, whenever $\eta_1\neq \eta_2\in S$ agree on $S^*$, they
have the same transition probabilities to any $\eta'\in S$, thus contributing a $0$-eigenvalue to $\bar{P}$. It is then
easy to see that all other eigenvalues of $\bar{P}$ (beyond those that originated from $\bar{P}^*$) must be 0. Altogether,
\begin{equation}\label{eq-barp*} \gap(\bar{P}^*) = \gap(\bar{P})\quad(\;=\bar{\gap}^\tau\;)~,\end{equation}
and it remains to give a lower bound for $\gap(\bar{P}^*)$. The next lemma shows that $\bar{P}^*$ is contracting with respect to
Hamming distance on $\{\pm1\}^{S^*}$.
\begin{lemma}\label{lem-contraction}
Let $\bar{X}^*_t$ and $\bar{Y}^*_t$ be instances of the chain $\bar{P}^*$, starting from $\varphi$ and $\psi$ respectively. Then there exists a coupling such that
$$\E_{\varphi, \psi} \dist (\bar{X}^*_1, \bar{Y}^*_1)\leq \bigg(\frac{b^{\lup}}{b^{\lup}+1} + \frac{1}{b^{\lup}+1} \cdot\frac{1 + (b-1)\lup }{b \kappa (1-\theta) (\ldown - \lup)}\bigg)\dist(\varphi, \psi)~.$$
\end{lemma}
\begin{proof}
Clearly, if $\varphi=\psi$ the lemma trivially holds via the identity coupling. In order to understand the setting when $\varphi\neq\psi$, recall the definition of the chain $\bar{P}^*$, which has the following two possible types of moves $E_1$ and $E_2$:
\begin{enumerate}
\item With probability $1 - \frac{1}{b^{\lup} + 1}$, the block dynamics updates one of the non-distinguished blocks: denote this event by $E_1$. Since this operation does not affect the value of the spins in the subset $S$ (and in particular, in $S^*$), the projection chain $\bar{P}$ remains in place in this case (and so does $\bar{P}^*$).
\item With probability $\frac{1}{b^{\lup} + 1}$, the distinguished block is being updated: denote this event by $E_2$. By the discussion above, this is equivalent to the following. Let $\eta$ denote the current state of the chain $\bar{P}^*$. First, $T\setminus S$ is assigned values according to the stationary distribution with boundary $\eta$ and $\tau$. Then, the distinguished block $B(\rho, \ldown)$ is updated given all other spins in the tree, and the resulting value of $S$ (and hence also of $S^*$) is determined by the new state of the projection chain.
\end{enumerate}
By the triangle inequality, it suffices to consider the case of $\dist(\varphi, \psi) =1$. Suppose therefore that $\varphi$ and $\psi$ agree everywhere on $S^*$ except at some vertex $\varrho$, and that without loss of generality,
 $$\varphi(\varrho) = 1~,~ \psi(\varrho) =-1~.$$
Crucially, the above mentioned procedure for the event $E_2$ is precisely captured by the spatial-mixing properties that were studied in Section \ref{sec:spatial-mixing}. Namely, a spin of some site $v \in S^*$ is propagated down the tree $T_v$ (with boundary condition $\tau$), and then the new value of $S^*$ is reconstructed from level $\ldown+1$, the external boundary of $B(\rho,\ldown)$. We construct a monotone coupling that will accomplish the required contraction property.

First, when propagating the sites $v\in S^*$ with $v\neq \varrho$, we use the identity coupling (recall that $\varphi(v)=\psi(v)$ for all $v\neq\varrho$). Second, consider the process that the spin at $\varrho$ undergoes. For $\varphi$, a positive spin is propagated to $T_\varrho$ (with boundary condition $\tau$) and then reconstructed from level $\ldown+1$ in the tree $T$ (which corresponds to level $\ldown-\lup+1$ in the subtree $T_\varrho$), with an additional boundary condition from $T \setminus T_\varrho$ that translates into some external field. For $\psi$, a negative spin is propagated analogously, and notice that in its reconstruction, the exact same external field applies (as $T \setminus T_\varrho$ was guaranteed to be the same for $\varphi$ and $\psi$).

Therefore, applying Proposition \ref{prop-delta} on the tree $T_\varrho$ with respect to the subtree $\hat{T}=B(\varrho,\ldown-\lup+1)$, we can deduce that
\begin{equation}\label{eq-dist-X-Y}\E_{\varphi, \psi}\Big(\bar{X}^*_1(\varrho)-\bar{Y}^*_1(\varrho) \given E_2\Big) \leq \frac{\mathrm{cap}_2 (B(\varrho, \ldown - \lup+1))}{\kappa (1-\theta)}~,\end{equation}
where $\kappa>\frac1{100}$, and the resistances are assigned as $$R_{(u,v)} = (\tanh \beta_c)^{-2\dist(\varrho,v)}~.$$
We now turn to estimating the $L^2$-capacity, which is equivalent to the effective conductance between $\varrho$ and $\partial B(\varrho, \ldown - \lup+1)$. This will follow from the well-known Nash-Williams Criterion (cf., e.g., \cite{LP}). Here and in what follows,
$R_\mathrm{eff} \deq 1/C_\mathrm{eff}$ denotes the effective resistance.
\begin{lemma}[Nash-Williams Criterion \cite{NashWilliams}]\label{lem-nw}
If $\{\Pi_j\}_{j=1}^J$ is a sequence of pairwise disjoint cutsets in a network $G$ that separate a vertex $v$ from some set $A$, then
$$R_\mathrm{eff}(v \leftrightarrow A) \geq \sum_j \Big(\sum_{e \in \Pi_j} \frac1{R_e}\Big)^{-1}~.$$
\end{lemma}
In our case, $G$ is the $b$-ary tree $B(\varrho,\ldown-\lup+1)$, and it is natural to select its different levels as the cutsets $\Pi_j$.
It then follows that
\begin{equation}\label{eq-R-inequality}
R_\mathrm{eff} \big(\varrho \leftrightarrow \partial B(\varrho, \ldown - \lup+1)\big) \geq \sum_{k=1}^{\ldown - \lup+1} (b^k \theta^{2k})^{-1} = \ldown - \lup+1~,\end{equation}
where we used the fact that $\tanh \beta_c = \theta = 1/\sqrt{b}$. It therefore follows that
$$\mathrm{cap}_2 \left(B(\varrho, \ldown - \lup+1)\right) \leq \frac{1}{\ldown - \lup}~,$$
which, together with \eqref{eq-dist-X-Y}, implies that
\begin{equation}\label{eq-difference-varrho}\E_{\varphi, \psi}\Big(\bar{X}^*_1(\varrho)-\bar{Y}^*_1(\varrho) \given E_2\Big) \leq \frac{1}{ \kappa (1-\theta) (\ldown - \lup)}~.\end{equation}
Unfortunately, asides from controlling the probability that the spin at $\varrho$ will coalesce in $\varphi$ and $\psi$, we must also consider
the probability that $\varrho$ would remain different, and that this difference might be propagated to other vertices in $S^*$ (as part of the update of $B(\rho, \ldown)$). Assume therefore that the we updated the spin at $\varrho$ and indeed
$\bar{X}^*_1(\varrho)\neq \bar{Y}^*_1(\varrho)$, and next move on to updating the remaining vertices of $S^*$.
Since our propagation processes corresponding to $\bar{X}^*$ and $\bar{Y}^*$ gave every vertex in $T\setminus T_\varrho$ the same spin,
it follows that each vertex $v \in S^*$, $v\neq \varrho$, has the same external field in $\bar{X}^*$ and $\bar{Y}^*$, with the exception of the effect of the spin at $\varrho$.

We may therefore apply the next lemma of \cite{BKMP}, which guarantees that we can ignore this mentioned common external field when bounding the probability of propagating the difference in $\varrho$.
\begin{lemma}[\cite{BKMP}*{Lemma 4.1}]\label{lem-external-field}
Let $T$ be a finite tree and let $v \neq w$ be vertices in $T$. Let $\{J_e \geq 0 : e \in E(T)\}$ be the interactions on $T$, and let $\{ H(u)\in\R  : u \in V(T)\}$ be an external field on the vertices of $T$. We consider the following conditional Gibbs measures:
\begin{align*}
\mu^{+, H}&: \mbox{the Gibbs measure with external field $H$ conditioned on $\sigma(v) = 1$.}\\
\mu^{-, H}&: \mbox{the Gibbs measure with external field $H$ conditioned on $\sigma(v) = -1$.}
\end{align*}
Then $\mu^{+, H}(\sigma(w)) - \mu^{-, H}(\sigma(w))$ achieves its maximum at $H \equiv 0$.
\end{lemma}
In light of the discussion above, Lemma \ref{lem-external-field} gives that
\begin{align*}\E_{\varphi, \psi}&\Big(\frac12\sum_{v\in S^*} (\bar{X}^*_1(v) - \bar{Y}^*_1(v)) \given E_2\Big)\\
& \leq \E_{\varphi, \psi}\Big(\bar{X}^*_1(\varrho)-\bar{Y}^*_1(\varrho) \given E_2\Big)\Big(1 + \sum_{k=1}^{\lup - 1} \frac{b-1}{b}b^k \theta^{2k}\Big) \\
&\leq \frac{1 + (b-1)(\lup - 1)/b}{ \kappa (1-\theta) (\ldown - \lup)} = \frac{1 + (b-1)\lup }{ b\kappa (1-\theta) (\ldown - \lup)}~,\end{align*}
 where in the first inequality we used the propagation property of the Ising model on the tree (Claim \ref{claim-Ising-representation}), and in
the second one we used the fact that $\theta=\tanh(\beta_c) = 1/\sqrt{b}$, as well as the estimate in \eqref{eq-difference-varrho}.

We conclude that there exists a monotone coupling of $\bar{X}^*_t$ and $\bar{Y}^*_t$ with
 $$\E_{\varphi, \psi}\Big(\dist(\bar{X}^*_1, \bar{Y}^*_1) \given E_2\Big) \leq \frac{1 + (b-1)\lup}{ b \kappa (1-\theta) (\ldown - \lup)}~,$$
 which then directly gives that
 $$\E_{\varphi, \psi}\big(\dist(\bar{X}^*_1, \bar{Y}^*_1)\big) \leq \frac{b^{\lup}}{b^{\lup}+1} + \frac{1}{b^{\lup}+1} \cdot\frac{1 + (b-1)\lup}{ b \kappa (1-\theta) (\ldown - \lup)}~,$$
as required.
\end{proof}

The above contraction property will now readily infer the required bound for the spectral gap of $\bar{P}^*$ (and hence also for $\bar{\gap}^\tau$).
\begin{proof}[\emph{\textbf{Proof of Lemma \ref{lem-gap-proj}}}]
The following lemma of Chen \cite{Chen} relates the contraction of the chain with its spectral gap:
\begin{lemma}[\cite{Chen}]\label{lem-spectra-contraction}
Let $P$ be a transition kernel for a Markov chain on a metric space $\Omega$. Suppose there exists a constant $\iota $ such that for each $x, y \in \Omega$, there is a
coupling $(X_1, Y_1)$ of $P(x, \cdot)$ and $P(y, \cdot)$ satisfying
\begin{equation}\label{eq-contraction-condition}
\E_{x,y} (\dist(X_1, Y_1))\leq \iota \, \dist(x, y)~.
\end{equation}
Then the spectral gap of $P$ satisfies $\gap \geq 1- \iota$.
\end{lemma}
By Lemma \ref{lem-contraction}, the requirement \eqref{eq-contraction-condition} is satisfied with
\begin{align*}
\iota &= \frac{b^{\lup}}{b^{\lup}+1} + \frac{1}{b^{\lup}+1} \cdot \frac{1+(b-1)\lup}{b\kappa(1-\theta)(\ldown-\lup)} ~,
\end{align*}
and hence
\begin{align}\gap(\bar{P}^*)  \geq 1-\iota
&= \frac1{b^\lup+1}\bigg(1-  \frac{1+(b-1)\lup}{b\kappa(1-\theta)(\ldown-\lup)}\bigg) \nonumber\\
&\geq \frac1{b^\lup+1}\bigg(1- \frac{\alpha}{\kappa(1-\theta)(1-2\alpha)}\bigg)~,\label{eq-gap-projection-relax}
\end{align}
where in the last inequality we increased $1+(b-1)\lup$ into $b\lup$ to simplify the final expression.
This lower bound on $\gap(\bar{P}^*)$ translates via
\eqref{eq-barp*} into the desired lower bound on the spectral gap of the projection chain, $\bar{\gap}^\tau$.
\end{proof}

We are now ready to provide a lower bound on the spectral gap of the block dynamics, $\gap_\mathcal{B}^\tau$, and
thereby conclude
the proof of Theorem \ref{thm-gap-block-dynamics}.
By applying Theorem \ref{thm-gap-decomposition} to our decomposition of the block dynamics chain $P_\mathcal{B}^\tau$,
\begin{equation}
  \label{eq-gap-b-decompose}
\gap_\mathcal{B}^\tau \geq \frac{\bar{\gap}}{3}~\wedge~ \frac{\bar{\gap}\cdot \gap_{\min}}{3\gamma + \bar{\gap}}~,
\end{equation}
where
$$\gap_{\min} \deq \min_{\eta\in  \{\pm1\}^S} \gap_\eta^\tau~,~\quad\gamma \deq \max_{\eta\in \{\pm1\}^S} \max_{x\in \Omega_\eta} \sum_{y\in \Omega\setminus \Omega_\eta} P_\mathcal{B}^\tau(x, y) ~.$$
Lemma \ref{lem-gap-restr} gives that $\gap_{\min} \geq 1/(b^\lup+1)$, and clearly, as the spins in $S$ can only change if the distinguished block is updated, $\gamma \leq  1/(b^\lup+1)$. Combining these two inequalities, we obtain that
\begin{align}\frac{\bar{\gap}\cdot \gap_{\min}}{3\gamma + \bar{\gap}} &=
\frac{\gap_{\min}}{1+3\gamma/\bar{\gap}}\geq\frac{1}{(b^\lup+1)+3/\bar{\gap}}
\geq \frac14\Big(\frac1{b^\lup+1} \wedge \bar{\gap}\Big)
\label{eq-gap-bar-formula}
\end{align}
with room to spare. Together with \eqref{eq-gap-b-decompose}, this implies that
$$ \gap_\mathcal{B}^\tau \geq \frac1{4(b^\lup+1)} ~\wedge~ \frac14\bar{\gap}~, $$
and Lemma \ref{lem-gap-proj} now gives that
$$ \gap_\mathcal{B}^\tau \geq \frac1{4(b^\lup+1)}\bigg(1 - \frac{\alpha}{\kappa(1-\theta)(1-2\alpha)}\bigg)~,$$
as required. This concludes the proof of Theorem \ref{thm-gap-block-dynamics}, and completes the proof of the upper bound on the mixing time.
\qed

\begin{remark}\label{rem-M-300}
Throughout the proof of Theorem \ref{thm-1} we modified some of the constants (e.g., \eqref{eq-gap-projection-relax}, \eqref{eq-gap-bar-formula} etc.) in order to simplify the final expressions obtained. By doing the calculations (slightly) more carefully, one can obtain
an absolute constant of about $300$ for the upper bound in Theorem \ref{thm-1}.
\end{remark}

\section{Lower bounds on the mixing time and inverse-gap}\label{sec:lower-bound}
In this section, we prove Theorem \ref{thm-tmix-trel-lower-bound}, which provides lower bounds on the inverse-gap and mixing time of the critical Ising model on a $b$-ary tree with free boundary. Throughout this section, let $b \geq 2$ be fixed, and set $\theta = \tanh \beta_c = \frac{1}{\sqrt{b}}$.

\subsection{Lower bound on the inverse-gap}\label{sec:lower-bound-trel}
The required lower bound will be obtained by an application of the Dirichlet form \eqref{eq-dirichlet-form}, using a certain weighted sum of the spins as the corresponding test function.
\begin{proof}[\emph{\textbf{Proof of Theorem \ref{thm-tmix-trel-lower-bound}, inequality \eqref{eq-trel-lower-bound}}}]
Let $T$ be a $b$-ary tree, rooted at $\rho$, with $h$ levels (and $n = \sum_{k=0}^h b^k$ vertices). We will
show that
$$\gap^{-1}  \geq \frac{b-1}{6b} n h^2~.$$
For simplicity, we use the abbreviation $d(v) \deq \dist(\rho,v)$, and define
$$g(\sigma) \deq \sum_{v \in T} \theta^{d(v)} \sigma(v)~\mbox{ for $\sigma\in\Omega$}~.$$
By the Dirichlet form \eqref{eq-dirichlet-form}, and since $P(\sigma,\sigma') \leq \frac1n$ for any $\sigma,\sigma'\in\Omega$
in the discrete-time dynamics, we have that
\begin{align*}
\mathcal{E}(g) &= \frac{1}{2} \sum_{\sigma, \sigma'} [g(\sigma) - g(\sigma')]^2 \mu(\sigma)P(\sigma, \sigma') \\
&\leq \frac{1}{2} \max_\sigma \sum_{\sigma'} [g(\sigma) - g(\sigma')]^2 \mu(\sigma)P(\sigma, \sigma') \leq \frac12\sum_{k=0}^{h} \frac{b^k}{n} (2\theta^{k})^2\leq \frac{2h}{n}~.
\end{align*}
On the other hand, the variance of $g$ can be estimated as follows.
\begin{align*}
\var_\mu g &= \var_\mu \Big(\sum_{v \in T} \theta^{d(v)}\sigma(v)\Big)= \sum_{u, w \in T} \theta^{d(u) + d(w)} \Cov_\mu(\sigma(u), \sigma(w))\\
&= \sum_{u,v, w\in T} \theta^{d(u) + d(w)} \Cov_\mu(\sigma(u), \sigma(w)) \one_{\{u \wedge w = v\}}~,
\end{align*}
where the notation $(u \wedge w)$ denotes their most immediate common ancestor (i.e., their common ancestor $z$ with the largest $d(z)$).
Notice that for each $v\in T$, the number of $u,w$ that are of distance $i,j$ from $v$ respectively and have $v=u\wedge w$ is
precisely $b^i \cdot (b-1)b^{j-1}$, since determining $u$ immediately rules $b^{j-1}$ candidates for $w$.
Furthermore, by Claim \ref{claim-Ising-representation} we have
$$ \Cov_\mu(\sigma(u),\sigma(w)) = \theta^{d(u)+d(w) -2d(v)}~,$$
and so
\begin{align*}
\var_\mu g &=\sum_{u, v, w\in T} \theta^{d(u) + d(w)} \theta^{d(u) + d(w) - 2 d(v)} \one_{\{u \wedge w = v\}}\\
&=\sum_{k=0}^h b^k\sum_{i=0}^{h-k}\sum_{j=0}^{h-k}  b^i (b-1)b^{j-1} \theta^{2k + i+j} \theta^{i+j}\\
&=\frac{b-1}{b}\sum_{k=0}^{h} (h-k)^2 = \frac{b-1}{6b} h(h+1)(2h+1) \geq \frac{b-1}{3b}h^3~.
\end{align*}
Altogether, we can conclude that
\begin{equation}
  \label{eq-gap-bound-thm2}
\gap \leq \frac{\mathcal{E}(g) }{\var_\mu g } = \frac{6b}{b-1}\cdot \frac{1}{n h^2} ~,
\end{equation}
as required.
\end{proof}

\subsection{Lower bound on the mixing-time}\label{sec:lower-bound-tmix}
In order to obtain the required lower bound on the mixing time, we consider a ``speed-up'' version of the dynamics, namely
a custom block-dynamics comprising a mixture of singletons and large subtrees. We will show that, even for this faster
version of the dynamics, the mixing time has order at least $n\log^3 n$.

Let $T$ be a $b$-ary tree with $h$ levels (and $n = \sum_{k=0}^h b^k$ vertices).
Consider two integers $1 \leq \lup<\ldown \leq h$, to be specified later.
For every $v\in H_\lup$, select one of its descendants in $H_\ldown$ arbitrarily, and denote it by $w_v$. Write $W = \{w_v : v\in H_\lup\}$ as the set of all such vertices. Further define
$$ B_v \deq \left(T_v \setminus T_{w_v}\right) \cup \{w_v\}\qquad(\mbox{for each $v \in H_\lup$})~.$$
The \emph{speed-up dynamics}, $(X_t)$, is precisely the block-dynamics with respect to
$$\mathcal{B} = \{ B_v : v \in H_\lup \} \cup \bigcup_{u \notin W} \{ u \}~.$$
In other words, the transition rule of the speed-up dynamics is the following:
\begin{enumerate}[(i)]
  \item Select a vertex $u \in V(T)$ uniformly at random.
  \item If $u\not\in W$, update this site according to the usual rule of the Glauber dynamics.
  \item Otherwise, update $B_v$ given the rest of the spins, where $v\in H_\lup$ is the unique vertex with $u=w_v$.
\end{enumerate}

\begin{figure}
\centering \includegraphics[width=3.5in]{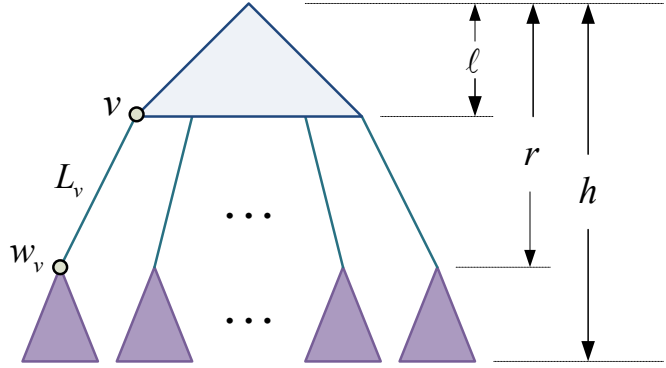}
\caption{Speed-up dynamics for the Ising model on the tree.}
\label{fig:speed-up-dynamics}
\end{figure}

The following theorem of \cite{PW} guarantees that, starting from all-plus configuration, the speed-up Glauber dynamics indeed mixes faster than the original one. In what follows, write $\mu \preceq \nu$ if $\mu$ stochastically dominates $\nu$. 
\begin{theorem}[\cite{PW} and also see \cite{Peres}*{Theorem 16.5}]\label{thm-PW-speed-up}
Let $(\Omega, S, V, \pi)$ be a monotone system and let $\mu$ be the distribution on $\Omega$ which results from successive updates at sites $v_1, \ldots, v_m$, beginning at the top configuration. Define $\nu$ similarly but with updates only at a subsequence $v_{i_1}, \ldots, v_{i_k}$. Then $\mu \preceq \nu$, and $\|\mu - \pi \|_{\mathrm{TV}} \leq \|\nu-\pi\|_{\mathrm{TV}}$. Moreover, this also holds if the sequence $v_1, \ldots, v_m$ and the subsequence $i_1, \ldots, i_k$ are chosen at random according to any prescribed distribution.
\end{theorem}
To see that indeed the speed-up dynamics $X_t$ is at least as fast as the usual dynamics, first note that any vertex $u\notin W$ is updated according to the original rule of the Glauber dynamics. Second, instead of updating the block $B_v$, we can simulate this operation by initially updating $w_v$ (given its neighbors), and then performing sufficiently many single-site updates in $B_v$. This approximates the speed-up dynamics arbitrarily well, and comprises a superset of the single-site updates of the usual dynamics. The above theorem thus completes this argument.

\begin{figure}
\centering \fbox{\includegraphics[width=4in]{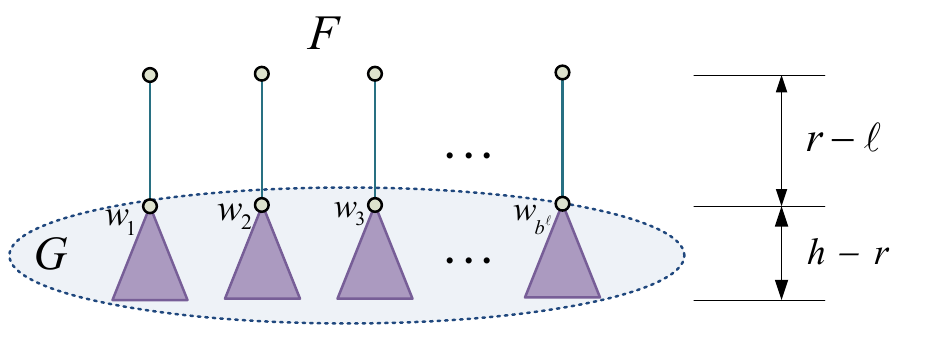}}
\caption{Speed-up dynamics on the forest $F$ and the sub-forest $G$.}
\label{fig:forest-dynamics}
\end{figure}

It remains to estimate the mixing time of the speed-up dynamics $X_t$. To this end, define another set of blocks as follows: for every $v\in H_\lup$, let $L_v$ denote the simple path between $v$ and $w_v$ (inclusive), define the forest
$$ F \deq \bigcup_{v \in H_\lup} \left(L_v \cup T_{w_v}\right)~,$$
and put
$$ \mathcal{B}_F \deq  \left\{ L_v : v\in H_\lup\right\} \cup \bigcup_{u \in F\setminus W} \{u\}~.$$
We define $Y_t$, the speed-up dynamics on $F$, to be the block-dynamics with respect to $\mathcal{B}_F$ above.
This should not be confused with running a dynamics on a subset of $T$ with a boundary condition of the remaining vertices;
rather than that, $Y_t$ should be thought of as a dynamics on a separate graph $F$, which is endowed with a natural
one-to-one mapping to the vertices of $T$.
Further note that, except for the singleton blocks in $\mathcal{B}$, every block $B_v\in\mathcal{B}$ in the block-dynamics
$X_t$ has a counterpart $L_v\subset B_v$ in $Y_t$.

The next lemma compares the continuous-time versions of $X_t$ and $Y_t$ (where each block is updated at rate 1),
and shows that on a certain subset of the vertices, they typically remain the same for a substantial amount of time.
\begin{lemma}\label{lem-speed-up-X-Y}
Let $(X_t)$ and $(Y_t)$ be the continuous-time speed-up dynamics on $T$ and $F$ respectively, as defined above.
Let $G = \bigcup_{v\in H_{\lup}}T_{w_v}$ and define
$$\tau= \inf_t\{X_t(u) \neq Y_t(u) \mbox{ for some } u\in V(G)\}~.$$
Then there exists a coupling of $X_t$ and $Y_t$ such that
$$\P(\tau > t) \geq  \exp(-\theta^{\ldown - \lup} b^{\lup}t)~.$$
\end{lemma}
\begin{proof}
For two configurations $\sigma\in \{\pm1\}^T$ and $\eta\in \{\pm1\}^F$, denote their Hamming distance on $F$ by
$$\dist(\sigma, \eta) = \sum_{v\in F} \one_{\{\sigma(v) \neq \eta(v)\}}~.$$
The coupling of $X_t$ and $Y_t$ up to time $\tau$ can be constructed as follows:
\begin{enumerate}
  \item Whenever a singleton block $\{u\}$ with $u\in T\setminus F$ is being updated in $X_t$, the chain $Y_t$ remains in place.
 \item Otherwise, when a block $B$ is updated in $X_t$, we update $B\cap F$ (the unique $B'\in\mathcal{B}_F$ with $B'\subset B$) in $Y_t$ so as to minimize $\dist(X_t, Y_t)$.
\end{enumerate}
For any $w\in W$, define the stopping time
$$\tau_w = \inf\{t: X_t(w) \neq Y_t(w)\} ~,$$
and notice that in the above defined coupling we have $\tau = \min_{w \in W} \tau_w$,
since $W$ separates $G\setminus W$ from $F$.

Let $v\in H_\lup$ and $w=w_v \in W$, and suppose that block $B_v$ is to be updated at time $t < \tau_w$ in $X_t$, and hence, as defined above, $L_v$ is to be updated in $Y_t$. By definition, at this time these two blocks have the same boundary except for at $v$, where there is a boundary condition in $T$ (the parent of $v$) and none in $F$ (recall $v$ is the root of one of the trees in $F$).

We now wish to give an upper bound on the probability that this update will result in $X_t(w) \neq Y_t(w)$.
By the monotonicity of the Ising model, it suffices to give an upper bound for this event in the case where $v$ has some parent $z$ in $F$, and $X_t(z) \neq Y_t(z)$.
In this case, we can bound the probability that $X_t(w) \neq Y_t(w)$ (in the maximal coupling) by an expression of the form
$$ \frac12\left(\mu^{+,H}(\sigma(w)) - \mu^{-, H}(\sigma(w))\right) $$ as described in Lemma \ref{lem-external-field}, where the external field $H$ corresponds to the value of the spins in $T_{w} \setminus \{w\}$. Lemma \ref{lem-external-field} then allows us to omit the
external field $H$ at $w$, translating
the problem into estimating the probability that a difference propagates from $v$ to $w$. By Claim \ref{claim-Ising-representation}, we deduce that
$$ \P \left(X_t(w) \neq Y_t(w)\right) \leq \theta^{\ldown-\lup}~,$$
and therefore
$$ \P \left(t < \tau_w\right) \geq \exp\left(- \theta^{\ldown-\lup}t\right)~.$$
Using the fact $|W| = b^\lup$, it follows that
$$\P(t < \tau) = \P\Big(t < \min_{w \in W}\tau_w\Big) \geq  \exp(-\theta^{\ldown - \lup} b^{\lup}t)~,$$
as required.
\end{proof}
With the above estimate on the probability that $X_t$ and $Y_t$ are equal on the subgraph $G$ up to a certain time-point,
we can now proceed to studying the projection of $X_t$ on $G$ via that of $Y_t$ (being a product chain, $Y_t$ is much simpler to analyze).

To be precise, let $\tilde{X}_t$ and $\tilde{Y}_t$ denote the respective projections of $X_t$ and $Y_t$ onto $G$,
which as a reminder is the union of all trees $T_{w_v}$. Notice that $\tilde{Y}_t$ is precisely the continuous-time single-site Glauber dynamics on $G$, since the block update of $L_v$ in $F$ translates simply into the single-site update of $w_v$ in $G$. On the other hand,
$\tilde{X}_t$ is not even necessarily a Markov chain. We next prove a lower bound on the mixing time of the Markov chain $\tilde{Y}_t$.

\begin{lemma}\label{lem-tilde-H-mu}
Let $\tilde{H}_t$ be the transition kernel of $\tilde{Y}_t$, and let $\mu_G$ denote its corresponding stationary measure.
Let $\gap'$ denote the spectral-gap of the continuous-time single-site dynamics on a $b$-ary tree of height $h-\ldown$.
Then
$$\|\tilde{H}_t (\allplus, \cdot) - \mu_G\|_\mathrm{TV} > \frac{3}{5}~\mbox{ for any $t \leq \frac{\lup \log b- 2}{2\gap'}$}~,$$
where $\allplus$ denotes the all-plus configuration.
\end{lemma}
\begin{proof}
Let $T'$ denote a $b$-ary tree of height $h-\ldown$ and $n'$ vertices. Let $P'$ be the transition kernel of the corresponding discrete-time single-site Glauber dynamics on $T'$, let $H'_t$ be the transition kernel of the continuous-time version of this dynamics, and let $\mu'$ be their
corresponding stationary measure.

By definition of $G$ as a disjoint union of $b^{\lup}$ copies of $T'$, clearly $\tilde{Y}_t$ is a product of $b^{\lup}$ copies of identical and independent component chains on $T'$. We can therefore reduce the analysis of $\tilde{Y}_t$ into that of $H'_t$, where the second eigenvalue of of its discrete-time counterpart $P'$ plays a useful role.

The following lemma ensures that $P'$ has an increasing eigenfunction corresponding to its second largest eigenvalue $\lambda'$.
\begin{lemma}[\cite{Nacu}*{Lemma 3}]
The second eigenvalue of the discrete-time Glauber dynamics for the Ising model has an increasing eigenfunction.
\end{lemma}
Since the eigenspace of $\lambda'$ has an increasing eigenfunction, it also contains a monotone eigenfunction $f$ such that $|f(\allplus)| = \|f\|_\infty$. Therefore, the transition kernel of the continuous-time chain satisfies
\begin{align}
\left(H'_t f\right)(\allplus) &= \bigg(\sum_{k=0}^{\infty}\mathrm{e}^{-t n'} \frac{(tn')^k}{k!}(P')^k f\bigg)(\allplus) \nonumber\\
&= \mathrm{e}^{- tn'}\sum_{k=0}^{\infty} \frac{(t n' \lambda')^k}{k!} f(\allplus) = \mathrm{e}^{-n' (1 - \lambda')t}f(\allplus)~.\label{eq-H'tf(1)}
\end{align}
Since $\int f d\mu' = 0$, we have that
\begin{align*}
|(H'_t f)(\allplus)| = \Big|\sum_y \left(H'_t (\allplus, y) f(y) - f(y)\mu'(y)\right)\Big|
&\leq 2\|f\|_\infty  \|H_t (\allplus, \cdot) - \mu'\|_{\mathrm{TV}}\,.
\end{align*}
Plugging in \eqref{eq-H'tf(1)} and using the fact that $|f(\allplus)| = \|f\|_\infty$, it follows that
\begin{equation}\label{eq-H'-TV}
\|H'_t (\allplus, \cdot) - \mu'\|_{\mathrm{TV}} \geq \frac{1}{2}\mathrm{e}^{-n'(1-\lambda')t}~.
\end{equation}

In order to relate the product chain $\tilde{Y}_t$ to its component chain $Y'_t$, we will consider the Hellinger distance between certain
distributions, defined next (for further details, see, e.g., \cite{LeCam}). First, define the \emph{Hellinger integral} (also known as the \emph{Hellinger affinity}) of two distribution $\mu$ and $\nu$ on $\Omega$ to be
$$I_\mathcal{H}(\mu, \nu) \deq \sum_{x\in \Omega} \sqrt{\mu(x) \nu(x)}~.$$
The \emph{Hellinger distance} is now defined as
$$d_\mathcal{H} (\mu, \nu) \deq \sqrt{2 - 2 I_\mathcal{H} (\mu, \nu)}~.$$
Clearly, for any two distributions $\mu$ and $\nu$,
\begin{equation*}
I_\mathcal{H}(\mu, \nu) = \sum_{x\in \Omega} \sqrt{\mu(x) \nu(x)} \geq \sum_{x \in \Omega} \mu(x) \wedge \nu(x) = 1 - \|\mu - \nu\|_{\mathrm{TV}}~,
\end{equation*}
and so $d_\mathcal{H}$ provides the following lower bound on the total variation distance:
\begin{equation}\label{eq-eq-Hellinger-affinity-bound-tv}\|\mu - \nu\|_{\mathrm{TV}} \geq
1 - I_\mathcal{H}(\mu,\nu) = \frac12 d^2_\mathcal{H}(\mu,\nu)~.\end{equation}
Furthermore, the Hellinger distance also provides an upper bound on $d_\mathrm{TV}$, as the next simple inequality (e.g., \cite{EKPS}*{Lemma 4.2 (i)}) shows:
\begin{equation}\label{eq-Hellinger-TV}
\|\mu - \nu\|_{\mathrm{TV}} \leq d_\mathcal{H} (\mu, \nu)~.
\end{equation}
To justify this choice of a distance when working with product chains, notice that
 any two product measures $\mu = \prod_{i=1}^n \mu^{(i)}$ and $\nu = \prod_{i=1}^n \nu^{(i)}$ satisfy
\begin{equation}\label{eq-product-Hellinger}
I_\mathcal{H}(\mu, \nu) = \prod_{i=1}^n I_\mathcal{H}(\mu^{(i)}, \nu^{(i)})~.
\end{equation}
Next, we consider the Hellinger integral of our component chains $H'_t$. Indeed, combining the definition of $d_\mathcal{H}$ with \eqref{eq-Hellinger-TV}, we get that
$$I_\mathcal{H}(H'_t(\allplus, \cdot), \mu') \leq 1- \frac{1}{2} \|H'_t(\allplus, \cdot) - \mu'\|_{\mathrm{TV}}^2 \leq  1 - \frac18\mathrm{e}^{-2 n'(1 - \lambda') t }~,$$
where the last inequality is by \eqref{eq-H'-TV}.
Therefore, applying \eqref{eq-product-Hellinger} to the product chain $\tilde{H}_t$ (the product of $b^{\lup}$ copies of $H'_t$), we can now deduce that
$$I_\mathcal{H}(\tilde{H}_t(\allplus, \cdot), \mu_G) \leq \bigg(1 - \frac18\mathrm{e}^{-2(1 - \lambda') t n'}\bigg)^{b^{\lup}}~.$$
At this point, \eqref{eq-eq-Hellinger-affinity-bound-tv} gives that
$$\|\tilde{H}_t(\allplus, \cdot) - \mu_G\|_{\mathrm{TV}} \geq 1 - \Big(1- \frac{\mathrm{e}^{-2(1 - \lambda')t n'}}{8}\Big)^{b^{\lup}}~.$$
Recall that by definition, $\gap'$ is the spectral-gap of $H'_t$, the continuous-time version of $P'$, and so $\gap' = n' (1-\lambda')$. Hence, if
$$t \leq \frac{\lup \log b- 2}{2\gap'}$$
then
$$\|\tilde{H}_t (\allplus, \cdot) - \mu_G\|_\mathrm{TV} \geq 1 - \exp\big(-\mathrm{e}^{2}/8\big) > \frac{3}{5}~,$$
as required.
\end{proof}

The final ingredient required is the comparison between $\mu_G$ (the Gibbs distribution on $G$),  and the projection of $\mu$ (the Gibbs distribution for $T$) onto the graph $G$. The following lemma provides an upper bound on the total-variation distance between these two measures.
\begin{lemma}\label{lem-hat-tilde-mu}
Let $\mu$ and $\mu_G$ be the Gibbs distributions for $T$ and $G$ resp., and let $\tilde{\mu}$ denote the projection of $\mu$ onto $G$, that is: $$\tilde{\mu} (\eta) = \mu (\{\sigma\in \{\pm1\}^T : \sigma_G = \eta\})\qquad \mbox{( for $\eta\in \{\pm1\}^G$ )}~.$$
Then $\|\mu_G - \tilde{\mu} \|_{\mathrm{TV}} \leq b^{2 \lup} \theta^{2(\ldown - \lup)}$.
\end{lemma}
\begin{proof}
Recalling that $G$ is a disjoint union of trees $\{ T_w : w \in W\}$, clearly the configurations of these trees are independent according to $\mu_G$. On the other hand, with respect to $\tilde{\mu}$, these configurations are correlated through their first (bottom-most) common ancestor. Further notice that, by definition, the distance between $w_i\neq w_j\in W$ in $T$ is at least $2(\ldown-\lup+1)$, as they belong
to subtrees of distinct vertices in $H_\lup$.

To bound the effect of the above mentioned correlation, we construct a coupling between $\mu_G$ and $\tilde{\mu}$ iteratively on the trees $\{ T_w : w \in W\}$, generating the corresponding configurations $\eta$ and $\tilde{\eta}$, as follows. Order $W$ arbitrarily as $W = \{w_1,\ldots,w_{b^\lup}\}$, and begin by coupling $\mu_G$ and $\tilde{\mu}$ on $T_{w_1}$ via the identity coupling. Now, given a coupling on $\cup_{i < k} T_{w_i}$, we extend the coupling to $T_{w_k}$ using a maximal coupling. Indeed, by essentially the same reasoning used for the coupling of the processes $X_t$ and $Y_t$ on $G$ in Lemma \ref{lem-speed-up-X-Y}, the probability that some already determined $w_i$ (for $i< k$) would affect $w_k$ is at most $\theta^{2(\ldown-\lup+1)}$. Summing these probabilities, we have that
$$\P\left(\eta_{T_{w_k}} \neq \tilde{\eta}_{T_{w_k}}\right) = \P\left(\eta(w_k) \neq \tilde{\eta}(w_k)\right) \leq (k - 1) \theta^{2(\ldown - \lup+1)}~.$$
Altogether, taking another union bound over all $k\in[b^\lup]$, we conclude that
$$\|\mu_G - \tilde{\mu}\|_\mathrm{TV} \leq \P(\eta \neq \tilde{\eta}) \leq b^{2 \lup} \theta^{2(\ldown - \lup)}~,$$
completing the proof.
\end{proof}
We are now ready to prove the required lower bound on $\tmix$.
\begin{proof}[\emph{\textbf{Proof of Theorem \ref{thm-tmix-trel-lower-bound}, inequality \eqref{eq-tmix-lower-bound}}}]
As we have argued above (see Theorem \ref{thm-PW-speed-up} and the explanation thereafter), it suffices to establish a lower bound on the mixing time of the speed-up dynamics $X_t$ on $T$. By considering the projection of this chain onto $G$, we have that
\begin{align*}
\|\P_\allplus(X_t \in \cdot) - \mu\|_{\mathrm{TV}} &\geq \|\P_\allplus(\tilde{X}_t \in \cdot)  - \tilde{\mu}\|_{\mathrm{TV}}~,
\end{align*}
and recalling the definition of $\tau$ as $\inf_t\{ (X_t)_G \neq (Y_t)_G\}$,
\begin{align*}
\|\P_\allplus(\tilde{X}_t \in \cdot)  - \tilde{\mu}\|_{\mathrm{TV}} &\geq  \|\P_\allplus(\tilde{Y}_t \in \cdot)  - \tilde{\mu}\|_{\mathrm{TV}} - \P(\tau \leq t)\\
&\geq \|\P_\allplus(\tilde{Y}_t \in \cdot)  - \mu_G\|_{\mathrm{TV}} - \P(\tau \leq t) - \|\mu_G - \tilde{\mu}\|_{\mathrm{TV}}~.
\end{align*}
Let $\gap$ and $\gap'$ denote the spectral-gaps of the continuous-time single-site dynamics on a $b$-ary tree
with $h$ levels and $h-r$ levels respectively (and free boundary condition), and choose $t$ such that
\begin{equation}\label{eq-t-selection}t \leq \frac{\lup \log b- 2}{2\gap'}~.\end{equation}
Applying Lemmas \ref{lem-speed-up-X-Y}, \ref{lem-tilde-H-mu} and \ref{lem-hat-tilde-mu}, we obtain that
\begin{align}\label{eq-p-allplus-lower-bound}\|\P_\allplus(X_t \in \cdot) - \mu\|_{\mathrm{TV}} &\geq \frac35 - \left(1-\exp(-\theta^{\ldown - \lup} b^{\lup}t)\right) - b^{2 \lup} \theta^{2(\ldown - \lup)}~.
\end{align}
Now, selecting $$\lup = \frac{h}5~\mbox{ and }~\ldown = \frac{4h}5~,$$ and recalling that $b\theta^2 =1$, we have that the last two terms in \eqref{eq-p-allplus-lower-bound} both tend to $0$ as $h\to\infty$, and so
$$\|P_\allplus(X_t \in \cdot) - \mu\|_{\mathrm{TV}} \geq \frac35 - o(1)~.$$
In particular, for a sufficiently large $h$, this distance is at least $1/\mathrm{e}$, hence by definition the continuous-time
dynamics satisfies $\tmix \geq t$. We can can now plug in our estimates for $\gap'$ to obtain the required lower bounds on $\tmix$.

First, recall that by \eqref{eq-gap-bound-thm2},
$$\gap' \leq \frac{6b}{b-1}\cdot\frac1{(h-\ldown)^2}~,$$
and so the following choice of $t$ satisfies \eqref{eq-t-selection}:
$$t \deq \frac{(b-1)}{12b}(h - \ldown)^2(\lup \log b- 2) ~.$$
It follows that the mixing-time of the continuous-time dynamics satisfies
$$\tmix \geq t \geq  \Big(\frac{(b-1)\log b}{1500\, b} + o(1)\Big) h^3 ~,$$
and the natural translation of this lower bound into the discrete-time version of the dynamics yields the lower bound in \eqref{eq-tmix-lower-bound}.

Second, let $g(h)$
be the continuous-time inverse-gap of the dynamics on the $b$-ary tree of height $h$ with free-boundary condition,
and recall that by Theorem~\ref{thm-1}, we have that $g$ is
polynomial in $h$. In particular,
$$ g(h) \leq C g(h/5)\quad\mbox{for some fixed $C > 0$ and all $h$.}$$
Since by definition $(\gap')^{-1} = g(h-\ldown)=g(h/5)$ and $\gap^{-1} = g(h)$, we can choose $t$
to be the right-hand-side of \eqref{eq-t-selection} and obtain that for any large $h$
$$\tmix \geq t \geq C' \gap^{-1} h\quad\mbox{for some $C' > 0$ fixed.}$$
Clearly, this statement also holds when both $\tmix$ and $\gap$ correspond to the discrete-time version of the dynamics, completing the proof.
%
%
\end{proof}

\section{Phase transition to polynomial mixing}\label{sec:phase-transition}
This section contains the proof of Theorem \ref{thm-2}, which addresses the near critical Ising model on the tree, and namely, the transition of its (continuous-time) inverse-gap and mixing-time from polynomial to exponential in the tree-height. Theorem \ref{thm-2} will follow directly from the next theorem:
\begin{theorem}\label{thm-near-critical}
Fix $b \geq 2$, let $\epsilon=\epsilon(h)$ satisfy $0 < \epsilon < \epsilon_0$ for a suitably small constant $\epsilon_0$,
and let $\beta = \arctanh\big(\sqrt{(1+\epsilon)/b}\big)$. The following holds for the continuous-time Glauber dynamics for the Ising model on the $b$-ary tree with $h$ levels at the inverse-temperature $\beta$:
\begin{enumerate}[(i)]
  \item For some $c_1 > 0$ fixed, the dynamics with free boundary satisfies
 \begin{align}
\gap^{-1} \geq c_1  \left((1/\epsilon)\,\wedge\, h\right)^2 (1 + \epsilon)^h~.  \label{eq-trel-transition-lower}
\end{align}
\item For some absolute constant $c_2 > 0$ and any boundary condition $\tau$
\begin{align}
\gap^{-1} \leq \tmix \leq  \mathrm{e}^{c_2 (\epsilon h + \log h)}~.  \label{eq-trel-transition-upper}
\end{align}
\end{enumerate}
\end{theorem}

Throughout this section, let $b\geq 2$ be some fixed integer, and let $T$ be a $b$-ary tree with height $h$ and $n$ vertices.
Define $\theta = \sqrt{(1+\epsilon)/b}$, where $\epsilon=\epsilon(n)$ satisfies $0 < \epsilon \leq \epsilon_0$ (for some
suitably small constant $\epsilon_0<\frac18$ to be later specified), and as usual
 write $\beta = \arctanh(\theta)$.

\begin{proof}[\emph{\textbf{Proof of Theorem \ref{thm-near-critical}}}]
The proof follows the same arguments of the proof of Theorems \ref{thm-1} and \ref{thm-tmix-trel-lower-bound}. Namely, the upper bound uses an inductive step using a similar block dynamics, and the decomposition of this chain to establish a bound on its gap (as in Section \ref{sec:thm-1}) via the spatial mixing properties of the Ising model on the tree (studied in Section \ref{sec:spatial-mixing}). The lower bound will again follow from the Dirichlet form, using a testing function analogous to the one used in Section \ref{sec:lower-bound}. As most of the arguments carry to the new regime of $\beta$ in a straightforward manner, we will only specify the main adjustments one needs to make in order to extend Theorems \ref{thm-1} and \ref{thm-tmix-trel-lower-bound} to obtain Theorem \ref{thm-near-critical}.

\subsection*{Upper bound on the inverse-gap}
Let $\frac1{100}< \kappa < 1$ be the universal constant that was introduced in Lemma \ref{lem-f-analysis} (and appears in Proposition \ref{prop-delta} and Theorem \ref{thm-log-likelihood-decay}), and define
\begin{equation*}
\epsilon_0 \deq \frac\kappa{20}\leq \frac18~.
\end{equation*}
As $b \geq 2$ and $\epsilon<\epsilon_0 \leq \frac18$, we have that $\theta \leq \frac{3}{4}$, hence Proposition~\ref{prop-delta}
and Theorem~\ref{thm-log-likelihood-decay} both hold in this supercritical setting. It therefore remains to extend the arguments in Section \ref{sec:thm-1} (that use Proposition \ref{prop-delta} as one of the ingredients in the proof of the upper bound on $\gap^{-1}$) to this new regime of $\beta$.

Begin by defining the same block dynamics as in \eqref{eq-block-dynamics-def}, only with respect to the following choice
of $\lup$ and $\ldown$ (replacing their definition \eqref{eq-lup-ldown-def}):
\begin{align}
  \alpha &\deq \epsilon_0 =\kappa/20~, \label{eq-def-alpha-2}\\
\lup &\deq \alpha \left[(1/\epsilon)\, \wedge\, h\right]~,~\ldown \deq h-\lup~.  \label{eq-def-lup-ldown-2}
\end{align}
Following the same notations of Section \ref{sec:thm-1}, we now need to revisit the arguments of Lemma \ref{lem-contraction}, and extend them to the new value of $\theta = \tanh \beta =\sqrt{(1+\epsilon)/b}$. This comprises the following two elements:
\begin{enumerate}
  \item Bounding the $L^2$-capacity $\mathrm{cap}_2 (B(\varrho, \ldown - \lup))$.
 \item Estimating the probability that a difference in one spin would propagate to other spins, when coupling
two instances of the chain $\bar{P}^*$.
\end{enumerate}

Recalling the Nash-Williams Criterion (Lemma \ref{lem-nw}) and its application in inequality \eqref{eq-R-inequality}, the effective resistance between $\varrho$ and $\partial B(\varrho, \ldown - \lup)$ is at least
$$ \sum_{k=1}^{\ldown-\lup+1} \big(b^k\theta^{2k}\big)^{-1} = \sum_{k=1}^{\ldown-\lup+1} (1 + \epsilon)^{-k} =  \frac1\epsilon \left(1-(1+\epsilon)^{-(\ldown-\lup+1)}\right)~,$$
which implies that
\begin{equation}
  \label{eq-cap2-supercritical-bound}
\mathrm{cap}_2\left(B(\varrho, \ldown - \lup+1 )\right) \leq \frac{\epsilon}{1-(1+\epsilon)^{-(\ldown - \lup)}}~.
\end{equation}
Now, if $\epsilon \geq 1/h$, we have
\begin{align*}
1-(1+\epsilon)^{-(\ldown - \lup)} &= 1-(1+\epsilon)^{-(h-2\alpha/\epsilon)} \geq
1-(1+\epsilon)^{-(1-2\alpha)/\epsilon} \\
&\geq 1-\mathrm{e}^{-(1-2\alpha)} \geq \frac{1-2\alpha}2~,
\end{align*}
where the last inequality uses the fact that $\exp(-x)\leq 1-x+\frac{x^2}2$ and that $\alpha>0$.
Similarly, if $\epsilon < 1/h$ then
\begin{align*}
1-(1+\epsilon)^{-(\ldown - \lup)} &= 1-(1+\epsilon)^{-h(1-2\alpha)} \geq
1-\mathrm{e}^{-\epsilon h (1-2\alpha)} \\
&\geq \epsilon h (1-2\alpha) - \frac{(\epsilon h(1-2\alpha))^2}2 \geq \epsilon h \frac{1-2\alpha}2~,
\end{align*}
where in the last inequality we plugged in the fact that $\epsilon h < 1$. Combining the last two
equations with \eqref{eq-cap2-supercritical-bound}, we deduce that
\begin{equation*}
\mathrm{cap}_2\left(B(\varrho, \ldown - \lup+1 )\right) \leq \frac{2\left(\epsilon \;\vee\; (1/h)\right)}{1-2\alpha}~.
\end{equation*}
Using \eqref{eq-dist-X-Y}, it then follows that
$$ \E_{\varphi, \psi}\Big(\bar{X}^*_1(\varrho)-\bar{Y}^*_1(\varrho) \given E_2\Big) \leq
\frac{2\left(\epsilon \;\vee\; (1/h)\right)}{\kappa (1 - \theta)(1 - 2\alpha)}~.$$
By repeating the next arguments of Lemma \ref{lem-contraction} (without any additional essential changes), we obtain that under the monotone coupling,
\begin{align*}
&\E_{\varphi, \psi}\Big(\dist(\bar{X}^*_1, \bar{Y}^*_1) \given E_2\Big) \leq \frac{2\left(\epsilon \;\vee\; (1/h)\right)}{ \kappa (1 - \theta)(1 - 2\alpha)}
\Big(1 + \sum_{k=1}^{\lup - 1} \frac{b-1}{b}b^k \theta^{2k}\Big) \\
 &= \frac{2\left(\epsilon \;\vee\; (1/h)\right)}{ \kappa (1 - \theta)(1 - 2\alpha)}
\Big(\frac1b + \frac{b-1}b \frac{(1+\epsilon)^{\alpha[(1/\epsilon)\, \wedge\, h]}-1}{\epsilon}\Big) \\
&\leq
  \frac{2\left(\epsilon\; \vee\; (1/h)\right)}{\kappa (1 - \theta)(1 - 2\alpha)} \frac{(1+\epsilon)^{\alpha [(1/\epsilon)\, \wedge\, h]}-1}{\epsilon}
\leq \frac{2\left(\epsilon\; \vee\; (1/h)\right)}{\kappa (1 - \theta)(1 - 2\alpha)} \frac{\mathrm{e}^{\alpha[1\,\wedge\, \epsilon h]}-1}{\epsilon}\\
 & \leq  \frac{2\left(\epsilon\; \vee\; (1/h)\right)}{\kappa (1 - \theta)(1 - 2\alpha)} \frac{2\alpha [1\,\wedge\, \epsilon h]}{\epsilon} = \frac{4\alpha}{\kappa (1 - \theta) (1 -
2\alpha)}~,
\end{align*}
where in the last line we used the fact that $\mathrm{e}^x-1 < 2x$ for all $0 \leq x \leq 1$.
Again defining $g_h = n_h \min_{\tau} \gap_h^\tau$, we note that all the remaining arguments in Section \ref{sec:thm-1} apply in our case without requiring any modifications, hence the following recursion holds for $g_h$:
\begin{equation}\label{eq-recursion-g-new}
g_h \geq c(\alpha) g_{\ldown} = c(\alpha) g_{h - \alpha[(1/\epsilon)\,\wedge\, h]}~,
\end{equation}
where $$c(\alpha) \deq \frac{1}{8}\Big(1 -\frac{4\alpha}{\kappa (1 -
\theta) (1 - 2\alpha)}\Big)~.$$
Recalling the definition \eqref{eq-def-alpha-2} of $\alpha$, since $\theta\leq\frac34$ and $\kappa < 1$ we have that
$$ \frac{4\alpha}{\kappa (1 - \theta) (1 - 2\alpha)} = \frac{2}{(1-\theta)(10-\kappa)} < \frac{8}{9}~,$$
and so $c(\alpha) > 0$. We now apply the next recursion over $g_{h_k}$:
$$h_0 = h ~,~ h_{k+1} = \begin{cases}
 h_k - (\alpha/\epsilon) &\mbox{ if }  h_{k} \geq (1/\epsilon)~,\\
 (1-\alpha) h_k &\mbox{ if }  h_{k} \leq (1/\epsilon)~.
\end{cases}$$
Notice that by our definition \eqref{eq-def-alpha-2}, we have $\epsilon < \epsilon_0 = \alpha$.
With this in mind, definition \eqref{eq-def-lup-ldown-2} now implies that for any $h > 1/\alpha$ we have $\lup,\ldown \geq 1$.
Thus, letting $K = \min \{k: h_k \leq 1/\alpha\}$, we can conclude from \eqref{eq-recursion-g-new} that
\begin{align*}
 g_{h_k} &\geq c(\alpha) g_{h_{k+1}} \mbox{ for all $k < K$, and hence }\\
g_h &\geq \left(c(\alpha)\right)^K g_{h_K}~.
\end{align*}
By the definitions of $h_k$ and $K$, clearly
$$K \leq \frac{\epsilon}{\alpha}h + \log_{1/(1-\alpha)}\left( h \,\wedge\,(1/\epsilon)\right) =  O(\epsilon h + \log h)~.$$
Since $h_K \leq 1/\alpha$, clearly $g_{h_K} > c'$ for some constant $c' = c'(\alpha) > 0$, giving
$$g_{h} \geq c' \left(c(\alpha)\right)^K \geq \mathrm{e}^{-M(\epsilon h + \log h)}$$
for some constant $M = M(\alpha) > 0$ and any sufficiently large $n$. By definition of $g_h$, this provides an upper bound
on $\gap^{-1}$, and as $\tmix= O\left(\gap^{-1} \log ^2 n \right)$ (see Corollary~\ref{cor-tmix-trel-log-2} in Section~\ref{sec:prelim}), we obtain the upper bound on $\tmix$ that appears in \eqref{eq-trel-transition-upper}.

\subsection*{Lower bound on the inverse-gap}
We now turn to establishing a lower bound on the inverse-gap. Define the test function $g$ to be
the same one given in Subsection \ref{sec:lower-bound-trel}:
$$g (\sigma) = \sum_{v\in T} \theta^{\dist(\rho,v)} \sigma(v)~.$$
By the same calculations as in the proof of Theorem \ref{thm-tmix-trel-lower-bound} (Subsection \ref{sec:lower-bound-trel}), we
have that \begin{equation}\mathcal{E}(g) \leq \frac12 \sum_{k=0}^h \frac{b^k}{n}(2\theta^k)^2 = \frac{2}{n} \sum_{k=0}^h (1+
\epsilon)^k = \frac{2}{n} \frac{(1+\epsilon)^{h+1}
-1}{\epsilon}~,\label{eq-E-g-epsilon}\end{equation} whereas
\begin{align}
\var_\mu (g) &=\frac{b-1}{b}\sum_{k=0}^{h} b^k \theta^{2k}
\Big(\sum_{i=0}^{h-k} b^i \theta^{2i}\Big)^2 =
\frac{b-1}{b}\sum_{k=0}^{h} (1+\epsilon)^k \Big(\sum_{i=0}^{h-k}
(1+\epsilon)^i\Big)^2\nonumber\\
&=\frac{b-1}{b}\sum_{k=0}^{h} (1+\epsilon)^k \Big(\frac{(1+\epsilon)^{h-k+1} -1}{\epsilon}\Big)^2\nonumber\\
&=\frac{b-1}{b\,\epsilon^2}\sum_{k=0}^{h} \Big((1+\epsilon)^{2h-k+2} - 2(1+\epsilon)^{h+1} +(1+\epsilon)^k\Big)\nonumber\\
 &=  \frac{b-1}{b\,\epsilon^3}\Big((1+\epsilon)^{2h+3} -
\left(2 h + 3\right)\epsilon(1+\epsilon)^{h+1} -1\Big)~.\label{eq-var-lower-bound-2}
\end{align}
When $\epsilon \geq 8/h$ we have
\begin{align*}
\frac12 (1+\epsilon)^{h+2} - (2h+3)\epsilon &\geq \frac\epsilon2(h+2)+\frac{\epsilon^2}2\binom{h+2}{2} - (2h+3)\epsilon \\
&\geq (h+2)\epsilon\left(\frac12 + \epsilon \frac{h+1}{4} - 2\right) \geq 4~,
\end{align*}
and therefore in this case \eqref{eq-var-lower-bound-2} gives
\begin{align}\label{eq-var-g-epsilon-large}
\var_\mu (g) &\geq \frac{b-1}{2b}\frac{(1+\epsilon)^{2h+3}}{\epsilon^3}~.
\end{align}
Combining \eqref{eq-E-g-epsilon} and
\eqref{eq-var-g-epsilon-large}, the Dirichlet form
\eqref{eq-dirichlet-form} now gives that
\begin{equation}
\gap \leq \frac{4b}{b-1}\frac{\epsilon^2}{n
(1 + \epsilon)^{h}}\quad\mbox{ for $ \epsilon \geq 8/h$}~.\label{eq-gap-epsilon-large}
 \end{equation}
On the other hand, when $0 \leq \epsilon < 8/h$ we still have $b \theta^2
\geq 1$ and hence
\begin{align*}
 \var_\mu (g)=
\frac{b-1}{b}\sum_{k=0}^{h} b^k \theta^{2k} \Big(\sum_{i=0}^{h-k}
b^i \theta^{2i}\Big)^2 \geq \frac{b-1}{b}\sum_{k=0}^{h}(h-k)^2 \geq
\frac{b-1}{3b} h^3.
\end{align*}
In addition, using the fact that the expression $[(1+\epsilon)^{h+1}-1]/\epsilon$ in \eqref{eq-E-g-epsilon}
is monotone increasing in $\epsilon$, in this case we have
$$\mathcal{E}(g) \leq \frac2n \frac{\left(1+(8/h)\right)^{h+1}-1}{8/h} \leq \mathrm{e}^7 h / n~,$$
where the last inequality holds for any $h \geq 20$.
Altogether, the Dirichlet form \eqref{eq-dirichlet-form} yields (for such values of $h$)
\begin{equation}\label{eq-gap-epsilon-small}
\gap \leq \frac{3\mathrm{e}^7 b}{b-1}\frac{1}{n h^2}\quad\mbox{ for $0< \epsilon \leq 8/h$}~.\end{equation}
Combining \eqref{eq-gap-epsilon-large} and
\eqref{eq-gap-epsilon-small}, we conclude that
$$\gap \leq  \frac{3\mathrm{e}^{15}b}{b-1}\left[n(1 + \epsilon)^{h} \left((1/\epsilon)\,\wedge\, h\right)^2\right]^{-1}~,$$
where we used the fact that $(1+\epsilon)^h \leq \mathrm{e}$.
This gives the lower bound on $\gap^{-1}$ that appears in \eqref{eq-trel-transition-lower}, completing the proof of
Theorem \ref{thm-near-critical}.
\end{proof}

\section{Concluding remarks and open problems}\label{sec:conclusion}
\begin{list}{\labelitemi}{\leftmargin=2em}
\item We have established that in the continuous-time Glauber dynamics for the critical Ising model on a regular tree with arbitrary boundary condition, both the inverse-gap and the mixing-time are polynomial in the tree-height $h$. This completes the picture for the phase-transition of the inverse-gap (bounded at high temperatures, polynomial at criticality and exponential at low temperatures), as conjectured by the physicists for lattices. Moreover, this provides the first proof of this phenomenon for any underlying geometry other than the complete graph.
\item In addition, we studied the near-critical behavior of the inverse-gap and mixing-time. Our results yield the critical exponent of $\beta-\beta_c$, as well as pinpoint the threshold at which these parameters cease to be polynomial in the height.
\item For further study, it would now be interesting to determine the precise power of $h$ in the order of each the parameters $\gap^{-1}$ and $\tmix$ at the critical temperature. In the free-boundary case, our lower bounds for these parameters in Theorem~\ref{thm-tmix-trel-lower-bound} provide candidates for these exponents:

\begin{question}
Fix $b\geq 2$ and let $\beta_c = \arctanh(1/\sqrt{b})$ be the critical inverse-temperature for the Ising model on a $b$-ary tree of height $h$. Does the corresponding continuous-time Glauber dynamics with free boundary condition satisfy $\gap^{-1} \asymp h^2$ and $\tmix \asymp h^3$?
\end{question}

\item Both at critical and at near-critical temperatures, our upper bounds for the inverse-gap and mixing-time under an arbitrary boundary condition matched the behavior in the free-boundary case. This suggests that a boundary condition can only accelerate the mixing of the dynamics, and is further supported by the behavior of the model
under the all-plus boundary, as established in \cite{MSW}. We therefore conjecture the following monotonicity of $\gap^{-1}$ and $\tmix$ with respect to the boundary condition:

\begin{conjecture}\label{conj-monotone}
Fix $b\geq 2$ and $\beta > 0$, and consider the Ising model on a $b$-ary tree with parameter $\beta$. Denote by $\gap$ and $\tmix$ the spectral-gap and mixing time for the Glauber dynamics with free boundary, and by $\gap^\tau$ and $\tmix^\tau$ those with boundary condition $\tau$. Then
\begin{equation*}
\gap \leq \gap^\tau~\mbox{ and }~\tmix \geq \tmix^\tau ~\mbox{ for any $\tau$}.
\end{equation*}
\end{conjecture}

\item A related statement was proved in \cite{Martinelli94} for two-dimensional lattices at low temperature: It was shown that, in that setting, the spectral-gap under the all-plus boundary condition is substantially larger than the spectral-gap under the free boundary condition.
In light of this, it would be interesting to verify whether the monotonicity property, described in Conjecture~\ref{conj-monotone}, holds for the Ising model on an arbitrary finite graph.

\end{list}



\begin{bibdiv}
\begin{biblist}

\bib{AF}{book}{
    AUTHOR = {Aldous, David},
    AUTHOR = {Fill, James Allen},
    TITLE =  {Reversible {M}arkov Chains and Random Walks on Graphs},
    note = {In preparation, \texttt{http://www.stat.berkeley.edu/\~{}aldous/RWG/book.html}},
}

\bib{Baxter}{book}{
   author={Baxter, Rodney J.},
   title={Exactly solved models in statistical mechanics},
   note={Reprint of the 1982 original},
   publisher={Academic Press Inc. [Harcourt Brace Jovanovich Publishers]},
   place={London},
   date={1989},
   pages={xii+486},
}

\bib{BKMP}{article}{
   author={Berger, Noam},
   author={Kenyon, Claire},
   author={Mossel, Elchanan},
   author={Peres, Yuval},
   title={Glauber dynamics on trees and hyperbolic graphs},
   journal={Probab. Theory Related Fields},
   volume={131},
   date={2005},
   number={3},
   pages={311--340},
}

\bib{BRZ}{article}{
   author={Bleher, P. M.},
   author={Ruiz, J.},
   author={Zagrebnov, V. A.},
   title={On the purity of the limiting Gibbs state for the Ising model on
   the Bethe lattice},
   journal={J. Statist. Phys.},
   volume={79},
   date={1995},
   number={1-2},
   pages={473--482},
}

\bib{CCCST1}{article}{
   author={Carlson, J. M.},
   author={Chayes, J. T.},
   author={Chayes, L.},
   author={Sethna, J. P.},
   author={Thouless, D. J.},
   title={Bethe lattice spin glass: the effects of a ferromagnetic bias and
   external fields. I. Bifurcation analysis},
   journal={J. Statist. Phys.},
   volume={61},
   date={1990},
   number={5-6},
   pages={987--1067},
}

\bib{CCCST2}{article}{
   author={Carlson, J. M.},
   author={Chayes, J. T.},
   author={Chayes, L.},
   author={Sethna, J. P.},
   author={Thouless, D. J.},
   title={Critical Behavior of the Bethe Lattice Spin Glass},
   journal={Europhys. Lett.},
   volume={106},
   date={1988},
   number={5},
   pages={355--360},
}

\bib{CCST}{article}{
   author={Carlson, J. M.},
   author={Chayes, J. T.},
   author={Sethna, J. P.},
   author={Thouless, D. J.},
   title={Bethe lattice spin glass: the effects of a ferromagnetic bias and
   external fields. II. Magnetized spin-glass phase and the de
   Almeida-Thouless line},
   journal={J. Statist. Phys.},
   volume={61},
   date={1990},
   number={5-6},
   pages={1069--1084},
}

\bib{Chen}{article}{
author= {Chen, Mu-Fa},
 year={1998},
 title= { Trilogy of couplings and general formulas for lower bound of spectral gap},
   conference={
   title={Probability towards 2000},
   address={New York},
      date={1995},
   },
book={
     series={Lecture Notes in Statist.},
     volume={128},
     publisher={Springer},
     place={New York},
   },
   date={1998},
   pages= {123--136},
}

\bib{DLP}{article}{
  title   = {The mixing time evolution of Glauber dynamics for the Mean-field Ising Model},
  author  = {Ding, Jian},
  author = {Lubetzky, Eyal},
  author = {Peres, Yuval},
  journal = {Comm. Math. Phys.},
  note = {to appear},
}

\bib{DS1}{article}{
   author={Diaconis, Persi},
   author={Saloff-Coste, Laurent},
   title={Comparison techniques for random walk on finite groups},
   journal={Ann. Probab.},
   volume={21},
   date={1993},
   number={4},
   pages={2131--2156},
}

\bib{DS2}{article}{
   author={Diaconis, Persi},
   author={Saloff-Coste, Laurent},
   title={Comparison theorems for reversible Markov chains},
   journal={Ann. Appl. Probab.},
   volume={3},
   date={1993},
   number={3},
   pages={696--730},
}

\bib{DS}{article}{
   author={Diaconis, P.},
   author={Saloff-Coste, L.},
   title={Logarithmic Sobolev inequalities for finite Markov chains},
   journal={Ann. Appl. Probab.},
   volume={6},
   date={1996},
   number={3},
   pages={695--750},
}

\bib{DS3}{article}{
   author={Diaconis, P.},
   author={Saloff-Coste, L.},
   title={Nash inequalities for finite Markov chains},
   journal={J. Theoret. Probab.},
   volume={9},
   date={1996},
   number={2},
   pages={459--510},
}

\bib{DL}{collection}{
   title={Phase transitions and critical phenomena. Vol. 20},
   editor={Domb, C.},
   editor={Lebowitz, J. L.},
   note={Cumulative author, title and subject index, including tables of
   contents, Vol.\ 1--19},
   publisher={Academic Press},
   place={San Diego, CA},
   date={2001},
   pages={vi+201},
}
		
\bib{EKPS}{article}{
   author={Evans, William},
   author={Kenyon, Claire},
   author={Peres, Yuval},
   author={Schulman, Leonard J.},
   title={Broadcasting on trees and the Ising model},
   journal={Ann. Appl. Probab.},
   volume={10},
   date={2000},
   number={2},
   pages={410--433},
}

\bib{HH}{article}{
    author={Hohenberg, P.C.},
    author={Halperin, B.I.},
    title={Theory of dynamic critical phenomena},
    year={1977},
    journal={Rev. Mod. Phys.},
    volume={49},
    number={3},
    pages={435--479},
}


\bib{Ioffe1}{article}{
   author={Ioffe, Dmitry},
   title={Extremality of the disordered state for the Ising model on general
   trees},
   language={English, with English and French summaries},
   conference={
      title={Trees},
      address={Versailles},
      date={1995},
   },
   book={
      series={Progr. Probab.},
      volume={40},
      publisher={Birkh\"auser},
      place={Basel},
   },
   date={1996},
   pages={3--14},
}

\bib{Ioffe2}{article}{
   author={Ioffe, Dmitry},
   title={On the extremality of the disordered state for the Ising model on
   the Bethe lattice},
   journal={Lett. Math. Phys.},
   volume={37},
   date={1996},
   number={2},
   pages={137--143},
}

\bib{JSTV}{article}{
   author={Jerrum, Mark},
   author={Son, Jung-Bae},
   author={Tetali, Prasad},
   author={Vigoda, Eric},
   title={Elementary bounds on Poincar\'e and log-Sobolev constants for
   decomposable Markov chains},
   journal={Ann. Appl. Probab.},
   volume={14},
   date={2004},
   number={4},
   pages={1741--1765},
}

\bib{LF}{article}{
    author={Lauritsen, Kent B{\ae}kgaard},
    author={Fogedby, Hans C.},
    title={Critical exponents from power spectra},
    journal={J. Statist. Phys.},
    volume={72},
    date={1993},
    number={1},
    pages={189--205},
}

\bib{LeCam}{book}{
   author={Le Cam, Lucien},
   title={Notes on asymptotic methods in statistical decision theory},
   publisher={Centre de Recherches Math\'ematiques, Universit\'e de
   Montr\'eal, Montreal, Que.},
   date={1974},
   pages={xvi+270},
}

\bib{LPW}{book}{
    author = {Levin, David A.},
    author = {Peres, Yuval},
    author = {Wilmer, Elizabeth},
    title =  {Markov Chains and Mixing Times},
    year  =  {2007},
    note = {In preparation, available at
    \texttt{http://www.uoregon.edu/\~{}dlevin/MARKOV/}},
}

\bib{Lyons}{article}{
   author={Lyons, Russell},
   title={The Ising model and percolation on trees and tree-like graphs},
   journal={Comm. Math. Phys.},
   volume={125},
   date={1989},
   number={2},
   pages={337--353},
}

\bib{LP}{book}{
author = {{R. Lyons with Y. Peres}},
title = {Probability on Trees and Networks},
publisher = {Cambridge University Press},
date = {2008},
note = {In preparation. Current version is available at \texttt{http://mypage.iu.edu/\~{}rdlyons/prbtree/book.pdf}},
}

\bib{Martinelli97}{article}{
   author={Martinelli, Fabio},
   title={Lectures on Glauber dynamics for discrete spin models},
   conference={
      title={Lectures on probability theory and statistics},
      address={Saint-Flour},
      date={1997},
   },
   book={
      series={Lecture Notes in Math.},
      volume={1717},
      publisher={Springer},
      place={Berlin},
   },
   date={1999},
   pages={93--191},
}

\bib{Martinelli94}{article}{
   author={Martinelli, F.},
   title={On the two-dimensional dynamical Ising model in the phase coexistence region},
   journal={J. Statist. Phys.},
   volume={76},
   date={1994},
   number={5-6},
   pages={1179--1246},
}

\bib{MSW}{article}{
   author={Martinelli, Fabio},
   author={Sinclair, Alistair},
   author={Weitz, Dror},
   title={Glauber dynamics on trees: boundary conditions and mixing time},
   journal={Comm. Math. Phys.},
   volume={250},
   date={2004},
   number={2},
   pages={301--334},
}


\bib{MY}{article}{
   author={Murakami, Atsushi},
   author={Yamasaki, Maretsugu},
   title={Nonlinear potentials on an infinite network},
   journal={Mem. Fac. Sci. Shimane Univ.},
   volume={26},
   date={1992},
   pages={15--28},
   issn={0387-9925},
}

\bib{Nacu}{article}{
  journal = {Probability Theory and Related Fields},
  volume  = {127},
  pages   = {177-185},
  year    = {2003},
  title   = {Glauber Dynamics on the Cycle is Monotone},
  author  = {Nacu, \c{S}.},
}

\bib{NashWilliams}{article}{
   author={Nash-Williams, C. St. J. A.},
   title={Random walk and electric currents in networks},
   journal={Proc. Cambridge Philos. Soc.},
   volume={55},
   date={1959},
   pages={181--194},
}

\bib{PP}{article}{
  author = {Pemantle, Robin},
  author = {Peres, Yuval},
  title = {The critical Ising model on trees, concave recursions and nonlinear  capacity},
  note = {preprint},
  year = {2005},
}

\bib{Peres}{article}{
    AUTHOR = {Peres, Yuval},
    conference = {
        title = {Lectures on ``Mixing for Markov Chains and Spin Systems''},
        address = {University of British Columbia},
        date = {August 2005},
        },
    note ={Summary available at \texttt{http://www.stat.berkeley.edu/\~{}peres/ubc.pdf}},
}

\bib{PW}{article}{
  author = {Peres, Yuval},
  author = {Winkler, Peter},
  title = {Can extra updates delay mixing?},
  note = {In preparation.},
}

\bib{Preston}{book}{
   author={Preston, Christopher J.},
   title={Gibbs states on countable sets},
   note={Cambridge Tracts in Mathematics, No. 68},
   publisher={Cambridge University Press},
   place={London},
   date={1974},
   pages={ix+128},
}

\bib{SaloffCoste}{article}{
   author={Saloff-Coste, Laurent},
   title={Lectures on finite Markov chains},
   conference={
      title={Lectures on probability theory and statistics},
      address={Saint-Flour},
      date={1996},
   },
   book={
      series={Lecture Notes in Math.},
      volume={1665},
      publisher={Springer},
      place={Berlin},
   },
   date={1997},
   pages={301--413},
}



\bib{Soardi93}{article}{
   author={Soardi, Paolo M.},
   title={Morphisms and currents in infinite nonlinear resistive networks},
   journal={Potential Anal.},
   volume={2},
   date={1993},
   number={4},
   pages={315--347},
   issn={0926-2601},
}

\bib{Soardi94}{book}{
   author={Soardi, Paolo M.},
   title={Potential theory on infinite networks},
   series={Lecture Notes in Mathematics},
   volume={1590},
   publisher={Springer-Verlag},
   place={Berlin},
   date={1994},
   pages={viii+187},
}

\bib{WH}{article}{
  title = {Universality in dynamic critical phenomena},
  author = {Wang, Fu-Gao}
  author = {Hu, Chin-Kun},
  journal = {Phys. Rev. E},
  volume = {56},
  number = {2},
  pages = {2310--2313},
  year = {1997},
}


\bib{Wilson}{article}{
   author={Wilson, David Bruce},
   title={Mixing times of Lozenge tiling and card shuffling Markov chains},
   journal={Ann. Appl. Probab.},
   volume={14},
   date={2004},
   number={1},
   pages={274--325},
}
\end{biblist}
\end{bibdiv}
\end{document}